\documentclass{amsart}
\usepackage{amssymb}
\usepackage{graphicx}
\usepackage{slashed}

\usepackage[colorlinks]{hyperref}
\newtheorem{thm}{Theorem}[section]
\newtheorem{cor}[thm]{Corollary}
\newtheorem{thm-def}[thm]{Theorem-Definition}
\newtheorem{lem}[thm]{Lemma}
\newtheorem{prop}[thm]{Proposition}
\theoremstyle{definition}

\theoremstyle{remark}
\newtheorem{rem}[thm]{Remark}
\numberwithin{equation}{section}

\newcommand{\tensor}{\otimes}

\newcommand{\vacB}{|0\rangle}
\newcommand{\vac}{\vacB}
\newcommand{\Cplx}{\mathbb C}
\newcommand{\CC}{\mathbb C}
\newcommand{\ZZ}{\mathbb Z}
\newcommand{\LL}{\mathbb L}
\newcommand{\VV}{\mathbb V}
\newcommand{\MM}{\mathbb M}

\newcommand{\WW}{\mathbb W}
\newcommand{\pole}{\Cplx}

\newcommand{\de}{\partial}

\newcommand{\cA}{\mathcal{A}}

\newcommand{\cD}{{\mathcal{D}}}

\newcommand{\cL}{{\mathcal{L}}}
\newcommand{\cM}{{\mathcal{M}}}
\newcommand{\cN}{\mathcal{N}}
\newcommand{\cO}{{\mathcal{O}}}

\newcommand{\cT}{{\mathcal{T}}}

\newcommand{\cW}{{\mathcal{W}}}

\newcommand{\cZ}{{\mathcal{Z}}}
\newcommand{\cE}{{\mathcal{E}}}
\newcommand{\cC}{{\mathcal{C}}}

\newcommand{\ops}[1]{{}_{({#1})} }

\newcommand{\id}{{\rm id\ }}

\newcommand{\fg}{\frak{g}}
\newcommand{\fa}{\frak{a}}
\newcommand{\fb}{\frak{b}}
\newcommand{\fm}{\frak{m}}
\newcommand{\fn}{\frak{n}}
\newcommand{\fh}{\frak{h}}
\newcommand{\fz}{\frak{z}}

\newcommand{\set}[1]{\left\{#1\right\}}
\newcommand{\Fields}{{Fields}}
\newcommand{\ghat}{{\hat{\fg}}}

\begin{document}

\pagestyle{plain}

\title{A chiral  Borel-Weil-Bott theorem}
% \title{ Algebras of twisted chiral differential operators and affine Lie algebras at the critical level }
\author{T.Arakawa \and F.Malikov}
\thanks{T.A.\ is partially  supported
by the JSPS Grant-in-Aid  for Scientific Research (B)
No.\ 20340007.}
\thanks{F.M. is partially supported by an NSF grant.}

\maketitle

\begin{abstract}
We compute the cohomology of modules over the algebra of twisted
chiral differential operators over the flag manifold. This is
applied to (1) finding the character of $G$-integrable irreducible
highest weight modules over the affine Lie algebra at the critical
level, and (2) computing a certain elliptic genus of the flag
manifold. The main tool is a result that interprets the
Drinfeld-Sokolov reduction as a derived functor.

\end{abstract}

\section{Introduction and the main result}
Let $G$ be a simple complex Lie group,
$\fg=\fn_+\oplus\fh\oplus\fn_-$ its Lie algebra, $X$ the
corresponding flag manifold. If $\lambda\in\fh^*$ is an integral
weight, denote by $\cL_\lambda$ the corresponding invertible
$G$-equivariant sheaf of $\cO_X$-modules and by $\cD^{\lambda}_X$
the algebra of twisted differential operators acting on
$\cL_\lambda$. The action of $\fg$ on $\cL_\lambda$ defines  a Lie
algebra morphism
\begin{equation}
\label{morph-g-dx-intro}
\fg\rightarrow \Gamma(X,\cD^{\lambda}_X)
\end{equation}
and  a localization functor \cite{BB1}
\begin{equation}
\label{bb-lok-funct-intr}
\Delta:\fg{\rm-Mod}\rightarrow\cD^{\lambda}_X{\rm-Mod},\;
\Delta(A)=\cD^{\lambda}_X\otimes_{\fg}A,
\end{equation}
which has proved of essence in representation theory and served as a template in modern
mathematical physics.

From various points of view, it is important to find a reasonable
analogue of functor (\ref{bb-lok-funct-intr}) in the case of the
affine Lie algebra, $\ghat$, a universal central extension of the
loop algebra $\fg\otimes\CC[t,t^{-1}]$. Various approaches have been
developed, see papers such as \cite{Kash,KashTan1,KashTan2,BD2,FG1},
all valuable in many respects yet deficient one way or another. We
would like to explore yet another proposal, which is based on
consistently replacing the notions of Lie or associative algebras
with that of a vertex algebra. This approach is  not a panacea
either, but it does lead to a pleasing result, and it gives answers
to a few natural questions arising independently of vertex algebras
(characters of irreducible $\ghat$-modules at the critical level) or
even representation theory (elliptic genera attached to flag
manifolds).

Constructed in \cite{MSV} is a sheaf of vertex algebras, $\cD^{ch}_X$,  known also   as an
algebra of {\it chiral differential operators} (CDO); this is an analogue of $\cD_X=\cD^0_X$
\footnote{This sheaf has been further studied in \cite{GMSII,AG}. The case of an arbitrary
underlying smooth variety is treated in \cite{BD1,GMS}}
Recently, a sheaf of {\it twisted chiral differential operators}, $\cD^{ch,tw}_X$, was proposed
\cite{AChM}; it is an analogue of not so much $\cD^{\lambda}_X$ as of its universal version,
where, roughly speaking, $\lambda$ becomes a variable. An analogue of (\ref{morph-g-dx-intro})
is a vertex algebra morphism
\begin{equation}
\label{vertex(g)-cd-chiral-intro}
\pi: V_{-h^{\vee}}(\fg)\rightarrow
\Gamma(X,\cD^{ch,tw}_X),
\end{equation}
where $V_{-h^{\vee}}(\fg)$ is the vertex algebra attached to $\ghat$ at the critical level
$-h^{\vee}$.

It is a peculiar property of the critical level that
$V_{-h^{\vee}}(\fg)$ acquires a big center
$\fz(V_{-h^{\vee}}(\fg))\subset V_{-h^{\vee}}(\fg)$, a fundamental
result of Feigin and Frenkel, \cite{FF2,F1,F3}. Analogously, the
center of the vertex algebra $\Gamma (X,\cD^{ch,tw}_X)$ equals $H_X$,
a commutative vertex algebra of differential polynomials on $\fh^*$.
Restricting morphism (\ref{vertex(g)-cd-chiral-intro}) to the center
one obtains
\begin{equation}
\label{center-to-center-intro}
\pi(\fz(V_{-h^{\vee}}(\fg))\subset
H_X.
\end{equation}

It is clear that a 1-dimensional representation of $H_X$ is the same
as a Laurent series $\nu(z)\in\fh^*((z))$. Let us introduce
$\cD^{ch,tw}_X{\rm-Mod_{\nu(z)}}$, the category of (sheaves of)
$\cD^{ch,tw}_X$-modules such that $H_X$ acts according to
$\nu(z)\in\fh^*((z))$. This is a reasonable analogue of
$\cD^{\lambda}_X{\rm-Mod}$.

For each $\nu(z)=\nu_0/z+\nu_{-1}+\nu_{-2}z+\cdots$, there is a functor \cite{AChM}
\begin{equation}
\label{zhu-intro} \cZ hu_{\nu(z)}: \cD^{\nu_0}_X{\rm-Mod}\rightarrow
\cD^{ch,tw}_X{\rm-Mod_{\nu(z)}}.
\end{equation}
In fact, one can prove \cite{AChM}, if some extra assumptions hold, that
$\cD^{ch,tw}_X{\rm-Mod_{\nu(z)}}$ has more than one, trivial, object if and only if $\nu(z)$ is
as demanded, and if so, then the functor $\cZ hu_{\nu(z)}$ is an equivalence of categories.

Now form the composition
\begin{equation}
\label{local-aff-intro} \cZ
hu_{\nu(z)}\circ\Delta:\fg{\rm-Mod}\rightarrow\cD^{ch,tw}_X{\rm-Mod}_{\nu(z)}.
\end{equation}
By virtue of
(\ref{vertex(g)-cd-chiral-intro},\ref{center-to-center-intro}), for
each $A\in \fg{\rm-Mod}$, $\cZ hu_{\nu(z)}\circ\Delta(A)$ is a sheaf
of $V_{-h^{\vee}}(\fg)$-modules, hence $\ghat$-modules, with central
character $\nu(z)\circ\pi$,  see (\ref{center-to-center-intro}),  and one can think of $\cZ
hu_{\nu(z)}\circ\Delta(A)$ as a localization of $\Gamma(X,\cZ
hu_{\nu(z)}\circ\Delta(A))$. This is why $\cZ
hu_{\nu(z)}\circ\Delta$ can be regarded as an affine version of
(\ref{bb-lok-funct-intr}).

Thus various $A\in \fg{\rm-Mod}$  serve to localize various $\ghat$-modules. For example, if
$\nu_0$ is dominant, and $M_{\nu_0}^c$ is the corresponding contragredient Verma module, then
$\cZ hu_{\nu(z)}\circ\Delta(M_{\nu_0}^c)$ is a localization of the Wakimoto module of critical
level with highest weight $\nu_0$, \cite{FF2,F1,F3,W}.

The most interesting example of such localization occurs when
$\nu_0$ is a regular dominant integral weight and $V_{\nu_0}$ is the
simple (finite dimensional) $\fg$-module with highest weight
$\nu_0$. In this case, $\cZ hu_{\nu(z)}\circ\Delta(V_{\nu_0})$ is a
$G$-equivariant sheaf, and so $\Gamma(X,\cZ
hu_{\nu(z)}\circ\Delta(V_{\nu_0}))$ is an object of $\ghat-{\rm
Mod}^{G[[t]]}_{\nu(z)\circ\pi}$, the category of those
$\ghat$-modules at the critical level with central character
$\nu(z)\circ\pi$, where the action can be integrated to an action of
$G[[t]]$.

In the impressive series of papers \cite{FG1,FG2,FG3}, Frenkel and Gaitsgory prove that
$\ghat-{\rm Mod}^{G[[t]]}_{\chi(z)\circ\pi}$ is a semi-simple category with a unique simple
object, $\VV_{\nu(z)}$, the  Weyl module with highest weight $\nu_0$ quotiented out by the
central character $\nu(z)\pi$. It follows that the cohomology $H^i(X,\cZ
hu_{\nu(z)}\circ\Delta(V_{\nu_0}))$ is a direct sum  of a number of copies of $\VV_{\nu(z)}$.
Here is the main result of the paper.

\begin{thm}
\label{main-theorem-intro}
Denote by $\cL_{\nu(z)}^{ch}$ the sheaf $\cZ
hu_{\nu(z)}\circ\Delta(V_{\nu_0})$ and let
\[
\chi(\cL_{\nu(z)}^{ch})=\sum_{i=0}^{{\rm dim}X}(-1)^i{\rm ch}
H^i(X,\cL_{\nu(z)}^{ch}),
\]
where ${\rm ch}$ stands for the formal character, cf.(\ref{what - the-char-form-is}).
Then
\[
\tag{1} \chi(\cL_{\nu(z)}^{ch})=\sum\limits_{w\in
W}(-1)^{\ell(w)}e^{w\circ \nu_0}\times \prod_{\alpha \in
       \widehat{\Delta}^{re}_+}
(1-e^{-\alpha})^{-1}.
\]
\[
\tag{2} H^i(X,\cL_{\nu(z)}^{ch})=\oplus_{w\in
W,\,l(w)=i}\VV_{\nu(z)}[\langle\nu_0-
w\circ\nu_0,\rho^{\vee}\rangle],
\]
where $W$ is the Weyl group of $\fg$, $l(w)$ is the length of $w\in W$, and $\VV_{\nu(z)}[m]$
stands for $\VV_{\nu(z)}$ as a $\ghat$-module with conformal filtration shifted by $m$,
cf.sect.~\ref{Gradings and character formulas}.
\end{thm}

This result is an extension (perhaps one can say ``chiralization'') of the Borel-Weil-Bott theorem,
\[
H^i(X,\cL_{\nu_0})=\left\{\begin{array}{ll}
V_{\nu_0}&\text{ if }i=0\\
0&\text{ otherwise},
\end{array}
\right.
\]
where $\cL_{\nu_0}=\Delta(V_{\nu_0})$ is the $G$-equivariant line bundle attached to the highest weight $\nu_0$;
note the appearance of the higher cohomology in our situation.

  Theorem~\ref{main-theorem-intro} (1) considerably
simplifies in the limit when $e^\alpha\mapsto 1$,
$e^{-\delta}\mapsto q$ (homogeneous grading):
\begin{equation}
\label{homo=vers-of=thm-1.1.2}  \chi(\cL_{\nu(z)}^{ch},q)={\rm
dim}V_{\nu_0}\prod_{j=1}^{+\infty}(1-q^j)^{-{\rm 2dim}X}.
\end{equation}

The assertions of Theorem~\ref{main-theorem-intro} can be interpreted as solutions of problems
stated independently of chiral differential operator algebra theory and interesting in their
own right.

 One such interpretation is the following character formula.
\begin{cor}
\label{char-form-cor-intro}
\[
{\rm ch}\VV_{\nu(z)}=\frac{\sum\limits_{w\in
W}(-1)^{\ell(w)}e^{w\circ \nu_0}}{ \prod_{\alpha\in
\Delta_+}(1-e^{-\langle \nu_0+\rho,\alpha^{\vee}\rangle\delta})
\prod_{\alpha\in
       \widehat{\Delta}^{re}_+}(1-e^{-\alpha})}.
\]
\end{cor}
Note  that the homogeneous grading specialization ($e^\alpha\mapsto
1$, $e^{-\delta}\mapsto q$) again makes this formula into an
infinite product
\begin{equation}
\label{homo=vers-char-form} {\rm dim}_q\VV_{\nu(z)}={\rm
dim}V_{\nu_0}\prod_{j=1}^{+\infty}(1-q^j)^{-{\rm
2dim}X}\prod_{\alpha\in\Delta_+}(1-q^{\langle\nu_0+\rho,\alpha^{\vee}\rangle})^{-1}.
\end{equation}
We would like to point out that this character formula is not new: in the case of $sl_2$ it was worked out in \cite{M}; in the full generality
it was first recorded  in \cite{A3}; and it immediately follows from Theorem 5 and formula (5.3) of \cite{FG2}. Our point is  not so much the formula
itself but the fact that it nicely fits in and follows from the proposed geometric framework.

To obtain another interpretation, introduce the following generating
function of locally free sheaves over $X$:
\begin{equation}
\label{sheaf-giving-ell-genus-intro}
\cE_\lambda=\cL_{\lambda}\otimes\left(\otimes_{n=1}^{\infty}\left(\oplus_{m=0}^{\infty}q^{nm}S^m\cT_{X}\right)\right)
\otimes\left(\otimes_{n=1}^{\infty}\left(\oplus_{m=0}^{\infty}q^{nm}S^m\Omega_{X}\right)\right)
\end{equation}
Formally expanding out we obtain
\[
\cE_\lambda=\cL_\lambda+q\cE_{\lambda,1}+q^2\cE_{\lambda,2}+\cdots
\]
Now define
\[
\chi(\cE_\lambda,q)=\chi(\cL_\lambda)+q\chi(\cE_{\lambda,1})+q^2\chi(\cE_{\lambda,2})+\cdots\in\ZZ[[q]],
\]
where $\chi(\cE_{\lambda,n})=\sum_i(-1)^{i}{\rm dim}H^i(X,\cE_{\lambda,n})$, the Euler
characteristic of $\cE_{\lambda,n}$.

It is easy to show, following \cite{BL}, that
\[
\chi(\cL_{\nu(z)}^{ch},q)=\chi(\cE_{\nu_0},q).
\]
On the other hand, it is known, \cite{BL}, see also some explanations in \cite{GM}, that
$\chi(\cL_{\nu(z)}^{ch},q)$ is a version of elliptic genus of $X$. More precisely,
an elliptic genus $g_Q(X,q)$ is attached \cite{HBJ} to a formal power series in $x$, $Q(X)$, that may also
depend, as it does in our situation, on $q$. We have (letting for simplicity $\nu(z)=0$)
\[
\chi(\cL_{\nu(z)=0}^{ch},q)=g_Q(X,q),\;Q(x)=\frac{x}{1-e^{-x}}\prod_{n=1}^\infty
 (1-q^ne^{-x})^{-1}(1-q^ne^{x})^{-1}.
\]

An alternative physics interpretation is obtained by recalling that
Witten \cite{Witt} has identified the cohomology vertex algebra
$\sum_i H^i(X,\cL_{\nu(z)=0}^{ch})$ with the chiral algebra of an
appropriate $(0,2)$-supersymmetric sigma-model on $X$. It follows
that $\chi(\cL_{0/z}^{ch},q)$ is the index of the corresponding
Dirac operator on the loop space $\cL X$:
\[
\chi(\cL_{0/z}^{ch},q)=\text{Ind}(\slashed{D},\cL X).
\]

 Theorem~\ref{main-theorem-intro} (1), or rather its corollary (\ref{homo=vers-of=thm-1.1.2}),
  is then a computation of
 either of these three, defined differently if at all but equal to each other,
 quantities.

  Note also that all coefficients of the genus computed in
 (\ref{homo=vers-of=thm-1.1.2})
 happen to be  positive. This is a bit mysterious; more examples of this positivity
 phenomenon can be found in \cite{GM}.

 As an aside, we would like to mention that the problem of
 computing
 the cohomology groups $H^i(X,\cL_{\nu(z)}^{ch})$, $0\leq i\leq{\rm
 dim}X$, has been around since \cite{MSV}, where
 what is denoted here by $\cL_{0/z}^{ch}$ was introduced and the cohomology found in the case of
 $\fg=sl_2$.
  The cohomological dimension zero case of the problem
 was worked out in \cite{AChM}. Also proved in \cite{AChM} is an
 extension of Theorem~\ref{main-theorem-intro} to not
 necessarily dominant highest weights $\nu_0$ in the case where
 $\fg=sl_2$; this is based on the earlier work \cite{M}.

 Our proof of Theorem~\ref{main-theorem-intro} involves the study of
 the Drinfeld-Sokolov reduction on an appropriate category $\cO$ at
 the critical level. The result we obtain may be of interest in its
 own right. The Drinfeld-Sokolov reduction is a version of
 semi-infinite cohomology. The latter has been known since its
 inception in \cite{Feig} to be a mixture of homology and
 cohomology; a refined treatment of this phenomenon can be found
 in \cite{V}. We find, somewhat unexpectedly, that the
 Drinfeld-Sokolov reduction, $H^{\infty/2+\bullet}_{DS}(L\fn_+,?)$,
%  at the critical level
is more like
 homology (cf.\ \cite{A2,A}):
 \begin{itemize}
 \item the functor $H^{\infty/2+i}_{DS}(L\fn_+,?)=0$ if $i>0$;
 \item
 the functor $H^{\infty/2+0}_{DS}(L\fn_+,?)$ is right
 exact, and the class of modules with Verma filtration is adapted to
 this functor;  and
 \item
  $H^{\infty/2-i}_{DS}(L\fn_+,?)$, $i>0$, is isomorphic to
 the derived functor $L^iH^{\infty/2+0}_{DS}(L\fn_+,?)$.
 \end{itemize}
 This is recorded in  the main body of the text as
 Theorem~\ref{dr-sok-as-derived}.
 
\bigskip

{\em Acknowledgments.} Part of the work was done when
T.A. was visiting  Academia Sinica, Taiwan, 
Weizmann Institute, Israel and 
Isaac Newton Institute, UK, and
when F.M. was visiting IPMU of Tokyo University and
IHES, France. We would like to thank these institutions for hospitality and excellent working conditions.
Our thanks also go to the referee of the paper for many valuable suggestions and, most importantly, for
locating an error in the original proof of Theorem~\ref{dr-sok-as-derived} (1).

\section{Vertex algebras and chiral differential operators:
examples} \label{Vertex algebras and chiral differential operators:
examples}

We will work over $\CC$; all vector spaces will actually be vector superspaces; if $V$ is a
vector space, $a\in V$, then by $\tilde{a}$ we shall denote the parity of $a$.

\subsection{Examples of vertex algebras}
\label{Examples of vertex algebras}
\subsubsection{Definition of a vertex algebra.}
\label{Definition of a vertex algebra.}

A {\em field} on a vector space $V$ is a formal series $$a(z) =
\sum_{n\in \ZZ} a_{(n)} z^{-n-1} \in ({\rm End} V)[[z, z^{-1}]]$$
such that for any $v\in V$ one has $a_{(n)}v = 0$ if $n\gg 0$.

Let $\Fields (V)$ denote the space of all fields on $V$.

A {\em vertex algebra}  is a vector space $V$ with the following
data:
\begin{itemize}
  \item a linear map $Y: V \to \Fields(V)$,  %  ({\rm End\  }V)[[z, z^{-1} ]]$,
    $V\ni a \mapsto a(z) = \sum_{n\in \ZZ} a_{(n)} z^{-n-1}$
  \item an even vector ${\bf 1}\in V$, called {\em vacuum vector}
  \item a linear operator $\de: V \to V$, called {\em translation operator}
\end{itemize}
that satisfy the following axioms:
\begin{enumerate}
  \item (Translation Covariance)

 $ (\de a)(z) = \de_z a(z)$

  \item (Vacuum)

  ${\bf 1}(z) = \id$;

  $a(z)\vac \in V[[z]]$ and $a_{(-1)}\vac = a$

  \item (Borcherds identity)
   \begin{align}
   \label{Borcherds-identity}
   &\sum\limits_{j\geq 0} {m \choose j} (a\ops{n+j} b )\ops{m+k-j}\\
   =& \sum\limits_{j\geq 0} (-1)^{j} {n \choose j}\{ a\ops{m+n-j} b\ops{k+j} - (-1)^{n+\tilde{a}\tilde{b}}
   b\ops{n+k-j} a\ops{m+j}
   \}\nonumber
   \end{align}
\end{enumerate}

% \end{defn}

A vertex algebra $V$ is {\em graded} if  $V = \oplus_{n\geq 0}V_n$
and for $a\in V_i$, $b\in V_j$ we have
$$a_{(k)}b \in V_{i+j - k -1}$$ for all $k\in \ZZ$. (We put $V_i  = 0$ for $i<0$.)

We say that a vector $v\in V_m$ has {\em conformal weight} $m$ and
write $\Delta_v = m$.

If $v\in V_m $ we denote $v_k  = v_{(k - m +1)}$, this is the
so-called conformal weight notation for operators. One has $$v_k V_m
\subset V_{m -k}.$$

 A {\em morphism} of vertex algebras is a map $f: V \to W$ that preserves vacuum and satisfies $f(v_{(n)}v') = f(v)_{(n)}f(v')$.

\subsubsection{Vertex algebra modules}
\label{Vertex algebra modules.}
 A {\em module} over a vertex algebra $V$ is a vector space $M$
together with a map
\begin{equation}
\label{def-vert-mod-1} Y^M: V \to \Fields(M),\;  a \to Y^M(a,z) =
\sum_{n\in \ZZ} a^M_{(n)}z^{-n-1},
\end{equation}
 that satisfy the following axioms:
\begin{enumerate}
  \item $\vac^M (z)   = \id_M  $
  \item (Borcherds identity)
  \begin{align}\label{def-vert-mod-2}
     &\sum\limits_{j\geq 0} {m \choose j} (a^{}_{(n+j)} b
     )^M_{(m+k-j)}\\
  = &\sum\limits_{j\geq 0} (-1)^{j} {n \choose j}\{ a^M_{(m+n-j)} b^M_{(k+j)} -
  (-1)^{n+\tilde{a}\tilde{b}}b^M_{(n+k-j)} a^M_{(m+j)} \}\nonumber
  \end{align}
\end{enumerate}

A module $M$ over a graded vertex algebra $V$ is called {\em graded}
if $M = \oplus_{n\geq 0} M_n$ with
 $v_{k}M_l  \subset M_{l-k}$  (assuming $M_{n} = 0$ for negative $n$) for all $v\in V$.
 Note that we have switched to the conformal weight notation.

An increasing filtration $\{F_nM,\;n\geq 0\}$, $\cup_{n}F_nM=M$, $F_nM=\{0\}$ if
$n\ll0$, is called {\em conformal} if
\begin{equation}
\label{def-of-conf-filtr}
v_{k}F_lM  \subset F_{l-k}M,\;\forall v\in V.
\end{equation}

Note that if $M$ is graded, then $\{F_nM=\bigoplus_{i=0}^{n}M_i\}$ is a conformal filtration, but we shall have a chance to encounter conformally filtered $V$-modules
that are not graded.

 A {\em morphism of modules} over a vertex algebra $V$ is a map $f: M \to N$
 that satisfies $f(v^M_{(n)}m) = v^N_{(n)}f(m)$ for $v\in V$, $m\in M$.

\subsubsection{Commutative vertex algebras.}
\label{Commutative_vertex_algebras} A vertex algebra is said to be
  {\em commutative} if $a_{(n)} b =0$ for $a$, $b$ in $V$ and $n\geq 0$.
 It is known that a commutative vertex algebra is the same
 as a commutative associative algebra with  derivation.

If $W$ is a vector space we denote by $H_W$ the algebra of
differential polynomials on $W$. As an associative algebra it is a
polynomial algebra in variables $x_i$, $\de x_i$, $\de^{(2)}x_i$,
$\dots$ where $\set{x_i}$ is a basis of $W^*$. A commutative vertex
algebra structure on $H_W$ is uniquely determined by attaching the
field
 $x(z) = e^{z\de}x_i$  to each $x_i$.

 $H_W$ is equipped with grading such that
 \begin{equation}
 \label{weights-0-1-diff-poly}
 (H_W)_0=\pole,\; (H_W)_1=W^*.
 \end{equation}

\subsubsection{Beta--gamma system or a CDO over an affine space.}
\label{beta-gamma-system} Let $U$ be a purely even vector space, $\{x_i\}\subset U^*$ and
$\{\partial_i\}\subset U$, $1\leq i\leq N$, a pair of dual bases. Denote by $\cD^{ch}(U)$ the
vertex algebra that is generated by the vector space $U\oplus U^*$ and relations
\begin{equation}
\label{def-beta-gamma-syst} x_{i(n)}x_j=\partial_{i(n)}\partial_j=\partial_{i(n+1)}x_j=0\text{
if }n\geq 0,\;
\partial_{i(0)}x_j=\delta_{ij}{\bf 1}.
\end{equation}
The Borcherds identity (\ref{Borcherds-identity}) implies the following commutation relations
\begin{equation}
\label{beta-gamma-commutat-erl} [x_{i(m)},x_{j(n)}]=[\partial_{i(m)},\partial_{j(n)}]=0,\;
[\partial_{i(m)},x_{j(n)}]=\delta_{ij}\delta_{m,-n+1}.
\end{equation}
This suggests an index shift $\partial_{in}=\partial_{i(n)}$, $x_{in}=x_{i(n-1)}$, which allows
to beautify the last relation as follows
\begin{equation}
\label{beta-gamma-commutat-erl-2} [\partial_{im},x_{jn}]=\delta_{ij}\delta_{m,-n}.
\end{equation}
As a vector space, $\cD^{ch}(U)$ is freely generated from ${\bf 1}$ by the family of pairwise
commuting  operators $\partial_{in-1}$, $x_{in}$ with $n\leq 0$. Thus
\begin{equation}
\label{poly-nature-of-beta-gamma} \cD^{ch}(U)\stackrel{\sim}{\rightarrow}\CC[\partial_{in-1},
x_{in};n\leq 0,1\leq i,j\leq N],\;{\bf 1}\stackrel{\sim}{\rightarrow}1\in\CC.
\end{equation}

This vertex algebra is  graded so that the degree of operators $\partial_{in}$, $x_{in}$ is
$(-n)$. In particular,
\begin{equation}
\label{weight-0-1-b-g} \cD^{ch}(U)_0=\pole[x_{10},...,x_{N0}],\;
\cD^{ch}(U)_1=\bigoplus_{j=1}^{N}(x_{j,-1}\cD^{ch}(\CC^n)_0\oplus
\partial_{j,-1}\cD^{ch}(\CC^n)_0).
\end{equation}
We tend to think of $x_{j0}{\bf 1}$ as  the function $x_j$ on $U$,
$\partial_{j,-1}{\bf 1}$ as the vector field $\partial/\partial
x_j$, $x_{j,-1}$ as the differential form $dx_j$ so that
$\cD^{ch}(U)_0$ is identified with functions on $U$ and
$\cD^{ch}(\CC^n)_1$ becomes $\cT_{U}(U)\oplus\Omega_U(U)$. We shall
soon make more sense out of this interpretation; in particular, we
shall see that $\cD^{ch}(U)$ is the space of global sections of a
sheaf of {\it chiral differential operators} (CDO) over $U$,
$\cD^{ch}_U$, and that the latter direct sum is a result of making
choices, but the exact sequence
\begin{equation}
\label{vert-algebr-ex-seq-first} 0\rightarrow\Omega_U(U)\rightarrow
\cD^{ch}(U)_1\rightarrow \cT_U(U)\rightarrow 0
\end{equation}
is natural; here $\Omega_U(U)\rightarrow \cD^{ch}(U)_1$ is defined
by
\[
f(x_1,...x_N)dx_j\mapsto x_{j,-1}f(x_{10},...x_{N0}){\bf 1},
\]
and $\cD^{ch}(U)_1\rightarrow \cT_U(U)$ is defined by
\[
x_{j,-1}f(x_{10},...x_{N0}){\bf 1}\mapsto 0, \partial_{j,-1}f(x_{10},...x_{N0}){\bf 1}\mapsto
f(x_{10},...x_{N0})\partial/\partial x_j.
\]

It is quite clear that the assignment $U\mapsto \cD^{ch}(U)$ defines a functor on the category
of affine spaces with affine isomorphisms for morphisms; in other words, a change of variables
$x_i\mapsto a_{is}x_s + b_i$ canonically lifts to an isomorphism of $\cD^{ch}(U)$.

\subsubsection{A super-version: Clifford algebra.}
\label{A super-version: Clifford algebra.}
The discussion in sect.~\ref{beta-gamma-system} is
easily carried over to the case where the purely even $U$ is replaced with a supervector space
of dimension $M|N$. We shall need the example of dimension $0|N$, so we define the Clifford
vertex algebra, $Cl(U)$, to be a vertex superalgebra that is generated by the {\it purely odd}
vector space $U\oplus U^*$ and relations
\begin{equation}
\label{def-b-c-syst} \phi^*_{(n)}\phi^*=\phi_{(n)}\phi=\phi_{(n+1)}\phi^*=0\text{ if }n\geq
0,\; \phi_{(0)}\phi^*=\langle\phi,\phi^*\rangle,\;\phi\in U, \phi^*\in U^*.
\end{equation}

Upon introducing $\phi_{n}=\phi_{(n)}$, $\phi^*_{n}=\phi^*_{(n-1)}$ the last relation becomes
\begin{equation}
\label{b-c-commutat-erl-2} [\phi_{m},\phi^*_{n}]=\langle\phi,\phi^*\rangle\delta_{m,-n}.
\end{equation}
As a vector space, $Cl(U)$ is freely generated from ${\bf 1}$ by the family of pairwise
commuting operators $\phi_{in-1}$, $\phi^*_{in}$ with $n\leq 0$; here $\{\phi_i\}$ and
$\{\phi^*_i\}$ is a pair of dual bases. Thus
\[
Cl(U)\stackrel{\sim}{\rightarrow}\CC[\phi_{in-1}, \phi^*_{in};n\leq
0,1\leq i,j\leq N]{\bf 1},
\]
where the polynomial ring is regarded as a superpolynomial ring, all generators being odd.

The functoriality of $U\mapsto Cl(U)$ is obvious.

\subsubsection{Affine vertex algebras.}
\label{Affine vertex algebras.} Let $\fg$ be a simple Lie algebra and $(.,.): S^2 \fg \to \CC$
the {\it normalized} invariant bilinear form on $\fg$, i.e., the form such that square of the
length of the longest root w.r.t the bilinear form induced on the dual Cartan subalgebra is 2.
The {\em affine Lie algebra} $\ghat$ associated with $\fg$ and $(.,.)$ is a central extension
of $\fg\tensor \CC[t,t^{-1}]+\CC K$ with bracket
\[
 [x\tensor t^n, y\tensor t^m ] = [x,y]
\tensor t^{m+n} + \delta_{n+m, 0} (x,y).
\]
Denote by $V_k(\fg)$ the vertex algebra generated by $\fg$ with relations
\begin{equation}
\label{def-af-vert} x_{(0)}y=[x,y],\; x_{(1)}y=k(x,y).
\end{equation}
Denote by $\ghat_{\geq}$ the subalgebra $\fg\otimes\CC[t]\oplus\CC k$ and by $\CC_k$ its
1-dimensional module, where $\fg\otimes\CC[t]$ acts as 0 and $K\mapsto k$. One has
\begin{equation}
\label{vert-aff-as-weyl} V_k(\fg)\stackrel{\sim}{\rightarrow}{\rm Ind}_{\ghat_{\geq}}^\fg\CC_k.
\end{equation}
The field attached to $x\in\fg$ is $x(z)=\sum_{n\in\ZZ}x_nz^{-n-1}$, where $x_n$ stands for
$x\otimes t^n\in\ghat$ and is regarded as an operator acting on ${\rm
Ind}_{\ghat_{\geq}}^\fg\CC_k$.

$V_k(\fg)$ is a graded vertex algebra, $V_k(\fg)=\oplus_{n\geq 0}V_k(\fg)_{n}$, the grading
being uniquely determined by the condition that the conformal weight of ${\bf 1}\stackrel{{\rm
def}}{=}1\in\CC_k$ be 0, and the corresponding degree of $x_n$ be $(-n)$.

One can likewise define a vertex algebra associated with any Lie
algebra $\fa$ and an invariant bilinear form on it, $(.,.)$. Since
in general there is no distinguished such form, not even up to
proportionality, we shall use the notation $V_{(.,.)}(\fa)$, or
 $V_0(\fa)$ if $(.,.)=0$.

\subsection{Wakimoto modules and algebras of chiral differential operators}

\label{Wakimoto modules and algebras of chiral differential operators}

Let $G$ be a simple  connected complex Lie group,
 $\fg$ the corresponding Lie algebra,
 $\fg=\fn_+\oplus\fh\oplus\fn_-$ a triangular decomposition,  $B$, $B_-$ resp., the
 subgroups of $G$ corresponding to $\fn_+\oplus\fh$, $\fn_-\oplus\fh$ resp., $U$, $U_-$
 the maximal unipotent subgroups of $B$ and $B_-$ resp. We will be interested in the flag
 manifold of $G$ to be denoted by $X$ and realized as $G/B_-$.

\subsubsection{Feigin-Frenkel-Wakimoto bozonization}
\label{Feigin-Frenkel-Wakimoto bozonization} The natural action $G\times G/B_-\rightarrow
G/B_-$ defines a Lie algebra homomorphism
\begin{equation}
\label{li-alg-act-flag}
 \fg\rightarrow \Gamma(X,\cT_X).
 \end{equation}
 Denote by $\Delta_+$ the set of positive (relative to the fixed triangular decomposition)
 roots of $\fg$. Let $w(\fn_+)\subset\fg$ be the maximal nilpotent subalgebra spanned by the
 root vectors with roots in $w(\Delta_+)$, $w$ being any element of the Weyl group $W$. The
 flag manifold has an atlas consisting of $U^w$-orbits, $U^w(w\overline{B_-})$,
 where $U^w$ is the maximal unipotent subgroup associated to $w(\fn_+)$. Each such orbit is
 a $U^w$-torsor, and in order to simplify the notation, we will identify  $U^w$ with $U^w(w\overline{B_-})$
 by sending ${\rm id}$ in the former to $w\overline{B_-}$ in the latter.

 Morphism (\ref{li-alg-act-flag}) defines morphisms
\begin{equation}
\label{li-alg-act-flag-cell-open}
 \fg\rightarrow \Gamma(U^w,\cT_X),\; w\in W.
 \end{equation}
Enter the vertex algebra $\cD^{ch}(U^w)$,  sect.~\ref{beta-gamma-system}. Since $V_k(\fg)$ is
generated by $\fg$, see sect.~\ref{Affine vertex algebras.}, it is natural to ask whether
(\ref{li-alg-act-flag-cell-open}) can be lifted to a vertex algebra morphism
$V_k(\fg)\rightarrow \cD^{ch}(U^w)$. Note that (\ref{vert-algebr-ex-seq-first}) implies a
lifting is determined by (\ref{li-alg-act-flag-cell-open}) modulo $\Gamma(U^w,\Omega_{U^w})$.

The answer is `yes' but only for $k=-h^{\vee}$, minus the dual Coxeter number. This is the
content of an important result of Feigin and Frenkel (building on earlier work of Wakimoto
\cite{W}.)

\begin{thm}(\cite{FF1,F3})
\label {exist=of-wak-mod} There is a unique lift of (\ref{li-alg-act-flag-cell-open}) to
$\fg\rightarrow \cD^{ch}(U^w)_1$ that extends to a vertex algebra morphism
\begin{equation}
\label{aff-vetr-alg-cdo-cell}\pi_w: V_{-h^{\vee}}(\fg)\rightarrow
\cD^{ch}(U^w).
\end{equation}
\end{thm}

The center of a vertex algebra $V$ is defined to be
\begin{equation}
\label{def-of-centre} \fz(V)=\{v\in V{\rm s.t. }v_{(n)}V=0, n\geq 0\}.
\end{equation}

It is a striking feature of the {\it critical level} $k=-h^{\vee}$ that at this level
$V_{-h^{\vee}}(\fg)$ acquires a big center, another important result of Feigin and Frenkel.

\begin{thm}(\cite{FF2,F3})
\label {exist=of-center-aff-crit}

(1) There are elements $p_i\in V_{-h^{\vee}}(\fg)_{d_i}$, $1\leq
i\leq {\rm rk}\fg$, such that $\fz(V_{-h^{\vee}}(\fg))$ is generated
by $\{p_1,...,p_r\}$.

(2) As a vertex algebra, $\fz(\fg)\stackrel{{\rm
def}}{=}\fz(V_{-h^{\vee}}(\fg))$ is isomorphic to the algebra of
differential polynomials on a ${\rm rk}\fg$-dimensional space.
\end{thm}

A useful complement to Theorem~\ref{exist=of-wak-mod}, also due to Feigin and Frenkel, is that
\begin{equation}
\label{in-wak-centre-zero} \pi_w(\fz(\fg))=0.
\end{equation}

\subsubsection{A CDO over the flag manifold.}
\label{A CDO over the flag manifold.} Since $U^w$ is an affine
space, the discussion in sect.~\ref{beta-gamma-system} allows us to
attach to each $U^w$ a graded vertex algebra $\cD^{ch}(U^w)=\bigoplus_{j\geq 0}\cD^{ch}(U^w)_j$. We shall now see how
these can be glued together into a sheaf over $X$.

First of all, we observe \cite{MSV} that there is a sheaf over $U^w$
of which $\cD^{ch}(U^w)$ is the space of global sections. For each
polynomial $f\in\CC[U^w]$, let $U^w_f$ be the Zariski
open subset of $U^w$ obtained by deleting the zero locus $\{f=0\}$. Define
\[
\cD^{ch}(U^w_f)=\CC[U^w]_{(f)}\otimes_{\CC[U^w]}\cD^{ch}(U^w),
\]
which makes sense due to (\ref{poly-nature-of-beta-gamma}), where we identifiy, as promised,
$\CC[x_{10},...,x_{N0}]$ with $\CC[x_{1},...,x_{N}]=\CC[U]$.

The fact of the matter is the remark, proved in \cite{MSV} and going
back to Feigin, that the vertex algebra structure on $\cD^{ch}(U^w)$
extends to that on $\cD^{ch}(U^w_f)$. It is clear that the assignment
$Y\mapsto \cD^{ch}(Y)$ is a sheaf of graded vertex algebras on $U^w$
such that its space of global sections is $\cD^{ch}(U^w)$. Denote
this sheaf $\cD^{ch}_{U^w}$; this is a CDO over $U^w$.

In this setting, (\ref{vert-algebr-ex-seq-first}) reads
\[
0\rightarrow\Omega_{U^w}\rightarrow
(\cD^{ch}_{U^w})_1\rightarrow\cT_{U^w}\rightarrow 0.
\]
In fact, by construction, this extension is split, but this
splitting is not natural, as we shall see in a second. In any case,
it defines a filtration,  $\Omega_{U^w}\subset (\cD^{ch}_{U^w})_1$,
and the corresponding $\operatorname{gr}(\cD^{ch}_{U^w})_1=\Omega_{U^w}\oplus
\cT_{U^w}$.

Now recall that $\{U^w\}$ is an atlas of $X$, and there are
transition functions
\[
\rho_{wv}: \Omega_{U^w}|_{U^w\cap
U^v}\stackrel{\sim}{\rightarrow}\Omega_{U^v}|_{U^w\cap U^v};\;
\cT_{U^w}|_{U^w\cap
U^v}\stackrel{\sim}{\rightarrow}\cT_{U^v}|_{U^w\cap U^v}.
\]
We would like to lift them to the CDOs $\cD^{ch}_{U^w}$, $w\in W$.
\begin{thm}( \cite{GMSII})
\label{exist-cdo-on-flag}

{\rm (1)} There exist
$\CC$-isomorphisms
\[
\hat{\rho}_{wv}: (\cD^{ch}_{U^w}|_{U^w\cap
U^v})_1\stackrel{\sim}{\rightarrow}(\cD^{ch}_{U^v}|_{U^w\cap U^v})_1
\]
that satisfy

{\rm (a)} they uniquely extend to  vertex algebra isomorphisms
\[
\hat{\rho}_{wv}: \cD^{ch}_{U^w}|_{U^w\cap
U^v}\stackrel{\sim}{\rightarrow}\cD^{ch}_{U^v}|_{U^w\cap U^v};
\]

{\rm (b)} they preserve the filtration and the corresponding graded
morphisms equal the classical $\rho_{wv}$;

{\rm (c)} the cocycle condition,
$\hat{\rho}_{wu}=\hat{\rho}_{vu}\circ\hat{\rho}_{wv}$, holds on
triple intersections $U^w\cap U^v\cap U^u$.

\flushleft{\rm (2)} Over $X$, there is a CDO, $\cD^{ch}_X$, such
that $\Gamma(U^w, \cD^{ch}_X)=\cD^{ch}(U^w)$ and
$\{\hat{\rho}_{wu}\}$ are  transition functions.

{\rm (3)} The vertex algebra morphisms (\ref{aff-vetr-alg-cdo-cell})
are compatible with $\{\hat{\rho}_{wu}\}$ and define a vertex
algebra morphism
\[
\pi: V_{-h^{\vee}}(\fg)\rightarrow \Gamma(X, \cD^{ch}_X),\;
\pi(\fz(V_{-h^{\vee}}(\fg))=0.
\]
{\rm(4)} The transition functions  $\{\hat{\rho}_{wv}\}$ are not unique but the equivalence class
of the corresponding CDO is independent of the choice.
\end{thm}
Note that the transition functions $\{\hat{\rho}_{wu}\}$ are not
$\cO_X$-linear, and so $\cD^{ch}_X$ is not a sheaf of
$\cO_X$-modules, but the filtration on $(\cD^{ch}_X)_1$ extends to
the entire $\cD^{ch}_X$ so that the corresponding graded object is a
sheaf of locally trivial $\cO_X$-modules. We have, cf.
(\ref{sheaf-giving-ell-genus-intro}),
\begin{equation}
\label{sheaf-giving-gr-obj-untw} {\rm
Gr}\cD^{ch}_X=\left(\otimes_{n=1}^{\infty}\left(\oplus_{m=0}^{\infty}q^{nm}S^m\cT_{X}\right)\right)\otimes
\left(\otimes_{n=1}^{\infty}\left(\oplus_{m=0}^{\infty}q^{nm}S^m\Omega_{X}\right)\right)
\end{equation}

\subsubsection{Wakimoto modules of critical
level.}\label{def-of-wak-of-crit-1} By pull-back, $\cD^{ch}(U^{{\rm
id}})$ is a $V_{-h^{\vee}}(\fg)$-module, hence a $\ghat$-module.
According to (\ref{in-wak-centre-zero} ), the center operates on
$\cD^{ch}(U^{{\rm id}})$ as 0; we shall call this module, following
\cite{FF1,F3}, the {\it Wakimoto module} of critical level and zero
central character and denote it by $\WW_{\nu(z)=0}$.

\subsubsection{A universal twisted CDO over the flag manifold.}
\label{A universal twisted CDO over the flag manifold.} The
constructions of sect.~\ref{A CDO over the flag manifold.} can be,
rougly speaking, deformed. For each integral weight $\lambda\in
P\in\fh^*$, denote by $\cL_{\lambda}$ the corresponding
$G$-equivariant invertible sheaf. The action of $G$ determines a Lie
algebra morphism
\begin{equation}
\label{fg-map-tdo}
 \fg\rightarrow \Gamma(X,\cD^{\lambda}_X),
\end{equation}
where $\cD^{\lambda}_X$ is the algebra of (twisted) differential operators acting on
$\cL_\lambda$. Trivializing $\cL_\lambda$ over $U^w$, the latter morphism becomes
\[
\fg\rightarrow \Gamma(U^w,\cT_X)\oplus \Gamma(U^w,\cO_X).
\]
It is clear that the choices can be made to ensure that this
morphism depepends on $\lambda$ polynomially, cf. a similar
discussion in \cite{BB2}, sect.~2.5. Thus we obtain a
morphism
\begin{equation}
\label{lie-alg-act-on-sheaf-univers} \fg\rightarrow
\Gamma(U^w,\cT_X)\oplus \Gamma(U^w,\cO_X)\otimes\CC[\fh^*].
\end{equation}
A vertex algebra version of this is as follows. Let $H_X$ be the
commutative vertex algebra of differential polynomials on $\fh^*$,
sect.~\ref{Commutative_vertex_algebras}; this is an analogue of
$\CC[\fh^*]$. Feigin and Frenkel proved \cite{FF1,F3} that
(\ref{lie-alg-act-on-sheaf-univers}) extends, for each $w$, to a
vertex algebra morphism
\begin{equation}
\label{aff-vetr-alg-cdo-cell-twist}\pi_w:
V_{-h^{\vee}}(\fg)\rightarrow \cD^{ch}(U^w)\otimes H_X \text{ s.t.
}\pi_w|_{\fz(V_{-h^{\vee}}(\fg))}:\fz(V_{-h^{\vee}}(\fg))\hookrightarrow
H_X,
\end{equation}
where $\cD^{ch}(U^w)\otimes H_X$ is the result of the well-known operation of tensor product of
vertex algebras, see e.g.\cite{K}. Note that the first assertion of
(\ref{aff-vetr-alg-cdo-cell-twist}) is a reasonably easy consequence of
(\ref{aff-vetr-alg-cdo-cell}), \cite{FBZ}.

In this case, too, the morphisms $\pi_w$ can be arranged into a
single morphism from $V_{-h^{\vee}}(\fg)$ to a certain sheaf of {\it
twisted chiral differential operators} (TCDO), $\cD^{ch,tw}_X$.

\begin{thm}(\cite{AChM})
\label{exist-tw-ch-do-flag} There is a sheaf of graded vertex algebras, $\cD^{ch,tw}_X$, over
$X$ such that

{\rm (1)} $\Gamma(U^w,\cD^{ch,tw}_X)$ is isomorphic to
$\cD^{ch}(U^w)\otimes H_X$;

{\rm (2)} the tautological embeddings $H_X\hookrightarrow
\cD^{ch}(U^w)\otimes H_X$ define an embedding $H_X\hookrightarrow
\cD^{ch,rw}_{X}$ as a constant subsheaf;  furthermore, this makes
$H_X$  the center of $\cD^{ch,tw}_X$;

{\rm (3)} $\cD^{ch}_X$ is isomorphic to $\cD^{ch,tw}_X$ modulo the
ideal generated by $H_X$;

{\rm (4)} the morphisms (\ref{aff-vetr-alg-cdo-cell-twist}) define a
graded vertex algebra morphism
\[
\pi:V_{-h^{\vee}}(\fg)\rightarrow \Gamma(X,\cD^{ch,tw}_X)\text{ s.t.
}\pi(\fz(V_{-h^{\vee}}(\fg))\subset H_X.
\]
\end{thm}

\subsubsection{Wakimoto modules.}\label{Wakimoto modules.}
By pull-back, each $\Gamma(U^w,\cD^{ch,tw}_X)$ is a
$V_{-h^{\vee}}(\fg)$-, hence $\ghat$-module. Call, following
\cite{FF1,F3}, $\Gamma(U^{\rm id},\cD^{ch,tw}_X)$ a {\it Wakimoto
module} with highest weight $(0,-h^{\vee})$ and denote it by
$\WW_{0,-h^{\vee}}$.

Note that $\WW_{0,-h^{\vee}}$ is different from a closely related
Wakimoto module {\it of critical level} and zero central character
$\WW_{\nu(z)=0}$ introduced in sect.~\ref{def-of-wak-of-crit-1}:
$\WW_{0,-h^{\vee}}$ are bigger than $\WW_{\nu(z)=0}$ because they
contain $H_X$ and, unlike $\WW_{\nu(z)=0}$, can be deformed away
from the critical level.

\subsubsection{Modules over the twisted CDO}\label{Modules over the twisted CDO.}
The twisted CDO $\cD^{ch,tw}_X$ is a deformation of $\cD^{ch}_X$
only morally: the vertex algebra axioms resist letting  $p\in H_X$
be a number. The situation changes pleasingly upon passing to
$\cD^{ch,tw}_X$-modules, where we shall at once get families of
modules depending on ${\rm rk}\fg$ functional parameters.

We will call a sheaf of vector spaces $\cM$ a $\cD^{ch,tw}_X$-{\em
module} if

(1) for each open $U\subset X$, $\Gamma(U,M)$ is a
$\Gamma(U,\cD^{ch,tw}_X)$-module;

(2) the restriction morphisms
$\Gamma(U,\cM)\rightarrow\Gamma(V,\cM)$, $V\subset U$, are
 are  $\Gamma(U,\cD^{ch,tw}_X)$-module morphisms, where the
 $\Gamma(U,\cD^{ch,tw}_X)$-module structure on $\Gamma(V,\cM)$ is that of the pull-back w.r.t. to
the restriction map
$\Gamma(U,\cD^{tw}_X)\rightarrow\Gamma(V,\cD^{tw}_X)$;

 (3) $\cM$ is conformally filtered, cf. (\ref{def-of-conf-filtr}), namely, there is an increasing sequence of
 subsheaves

\begin{equation}
\label{filtr-module-1} \{{\rm
F}_n\cM,n\in\ZZ\},\;\cup_{n=-\infty}^{+\infty}{\rm F}_n\cM_n=\cM,\,
{\rm F}_n\cM\subset {\rm F}_{n+1}\cM,\, {\rm F}_n\cM=\{0\}\text{ if
}n\ll0
\end{equation}
so that
\begin{equation}
 \label{filtr-module-2}
(\cD^{ch,tw}_X)_{l}\cM_n\subset\cM_{n-l}.
\end{equation}

\bigskip
Denote by $\cD^{ch,tw}_X-{\rm Mod}$ the category of
$\cD^{ch,tw}_X$-modules.

 Since the vertex
algebra $H_X$ is commutative, its irreducibles are all 1-dimensional
and are in 1-1 correspondence with the algebra of Laurent series
with values in $\fh^*$. Specifically, if $\nu(z)\in \fh^*((z))$,
then the character $\CC_{\nu(z)}$ is a 1-dimensional $H_X$-module
defined by
\begin{equation}
\label{def-of-char} \nu: H_X\rightarrow
Fields(\CC_{\nu(z)}),\;H_X\ni\lambda\mapsto\lambda^{\CC_{\nu(z)}}(z)\stackrel{{\rm
def}}{=}\lambda(\nu(z)),
\end{equation}
cf. sect.~\ref{Vertex algebra modules.} and recall that $H_X$ is the
algebra of differential polynomials on $\fh^*$. For example, if
$\lambda\in \fh$, thus $\lambda$ is a linear function on $\fh^*$,
and $\nu(z)=\sum_n\nu_{n}z^{-n-1}$, then
\[
\nu(\lambda)(z)=\sum_{n\in \ZZ}\lambda(\nu_n)z^{-n-1}\text{ or
}\nu(\lambda)_{(n)}=\lambda(\nu_n).
\]

Denote by $\cD^{ch,tw}_X-{\rm Mod_{\nu(z)}}$ the full subcategory of
$\cD^{ch,tw}_X-{\rm Mod}$ consisting of those
$\cD^{ch,tw}_X$-modules, where $H_X$ acts according to the character
$\nu(z)$.

We will say that a character $\nu(z)\in \fh^*((z))$ has  {\em
regular singularity} if $\nu(z)=\nu_0
z^{-1}+\nu_{-1}+\nu_{(-2)}z+\cdots$.

It is easy to see \cite{AChM} that if $\nu(z)$ has regular
singularity, then for each $\cM\in\cD^{ch,tw}_X-{\rm Mod_{\nu(z)}}$,
\[
{\rm Sing}\cM\stackrel{{\rm def}}{=}\{m\in\cM\text{ s.t.
}v_{n}m=0\,\forall v\in\cD^{ch,tw}_X, n>0\}
\]
is naturally a $\cD^{\nu_0}_X$-module, hence a functor
\begin{equation}
\label{forg-functor} {\rm Sing}:\;\cD^{ch,tw}_X-{\rm
Mod_{\nu(z)}}\rightarrow \cD^{\nu_0}_X-{\rm Mod},\; \cM\mapsto{\rm
Sing}\cM.
\end{equation}
This functor has the left adjoint\cite{AChM}
\begin{equation}
\label{left-adj-to-forg-functor} \cZ hu_{\nu(z)}:\cD^{\nu_0}_X-{\rm
Mod}\rightarrow \cD^{ch,tw}_X-{\rm Mod_{\nu(z)}},
\end{equation}
which, as the notation suggests, is closely related to Zhu's work
\cite{Zhu}. Its construction is simple enough: define
\begin{equation}
\label{constr-of-zhu} \cZ
hu_{\nu(z)}(\cA)(U^w)=\cD^{ch}(U^w)\otimes_{\CC}\Gamma(U^w,\cA)
\end{equation}
and then mimick the proof of Theorem~\ref{exist-tw-ch-do-flag} to
glue these pieces together. In particular, one sees that even though
$\cZ hu_{\nu(z)}(\cA)(U^w)=\cD^{ch}(U^w)\otimes_{\CC}\Gamma(U^w,\cA)$ appears
graded by setting $\cZ
hu_{\nu(z)}(\cA)(U^w)_n=\cD^{ch}(U^w)_n\otimes_{\CC}\Gamma(U^w,\cA)$, the
actual sheaf is only filtered by
\begin{equation}
\label{constr-of-zhu-filtr} \Gamma(U^w,\cZ
hu_{\nu(z)}(\cA)_n)=\oplus_{j=0}^{n}\cD^{ch}(U^w)_j\otimes_{\CC}\Gamma(U^w,\cA)
\end{equation}

\begin{thm}(\cite{AChM} Theorem 5.2, Remark 5.4.)
\label{equi-of-cat-from-prev}

The functors
\[
\cZ hu_{\nu(z)}:\cD^{\nu_0}_X-{\rm Mod}\leftrightarrow
\cD^{ch,tw}_X-{\rm Mod_{\nu(z)}}:{\rm Sing}
\]
are quasiinverses of each other that establish an equivalence of categories.
\end{thm}

Note that for each $\cM\in\cD^{ch,tw}_X-{\rm Mod_{\nu(z)}}$, the
corresponding graded object, ${\rm Gr}_{F}\cM$ is an object of the
category $\cM\in\cD^{ch,tw}_X-{\rm Mod_{\nu_0/z}}$. (Indeed, the
character $\nu(z)=\nu_0/z+\nu_{-1}+\cdots$ is the only source of
inhomogeneity, see e.g. (\ref{constr-of-zhu}), and the passage to
the graded object object replaces $\nu(z)$ with $\nu_0/z$, a
homogeneous character.)

Having refined the
filtration (\ref{constr-of-zhu-filtr}) further,  as in the discussion that led to  (\ref{sheaf-giving-gr-obj-untw}),  one obtains
\begin{equation}
\label{sheaf-giving-gr-obj-tw} {\rm gr}\cZ
hu_{\nu(z)}(\cA)=\cA\otimes\left(\otimes_{n=1}^{\infty}\left(\oplus_{m=0}^{\infty}q^{nm}S^m\cT_{X}\right)\right)\otimes
\left(\otimes_{n=1}^{\infty}\left(\oplus_{m=0}^{\infty}q^{nm}S^m\Omega_{X}\right)\right),
\end{equation}
where we have used ${\rm gr}\cZ hu(\cA)$ in place of a more logical
but awkward ${\rm Gr}{\rm Gr}_{F}\cA$.

We shall use Theorem~\ref{equi-of-cat-from-prev} only as a source of
examples of modules from $\cD^{ch,tw}_X-{\rm Mod_{\nu(z)}}$.

\subsubsection{Examples of $\cD^{ch,tw}_X$-modules.}
\label{examples-of-d-ch-tw-mod} To obtain examples of
$\cD^{ch,tw}_X$-modules with central character $\nu(z)$, we need a
supply of $\cD^{\nu_0}_X$-modules. For the purposes of
representation theory, the most interesting are those
$\cD^{\nu_0}_X$-modules that are $U$-equivariant, hence
supported on a union of  $U$-orbits.

Let $X_w\stackrel{{\rm def}}{=}Uw\overline{B_-}$, $w\in W$,  $i_w: X_w\hookrightarrow X$ the
tautological embedding, $i_{w,+}:\cD_{X_w}-{\rm Mod}\rightarrow \cD_{X}-{\rm Mod}$ the
$\cD$-module direct image functor. We obtain a family of $\cD_X$-modules, $i_{w,+}\cO_{X_w}$,
$w\in W$; $i_{w,+}\cO_{X_w}$ is often referred to as ``the module of distributions supported on
$X_w$.''

The space $\Gamma(U^w, i_{w,+}\cO_{X_w})$, which is essentially the space of global sections,
is easy to describe. Identify $U^w$ with $w(\fn_+)$ be means of the exponential map. Let
$\{x_\alpha,\,\alpha\in w(\Delta_+)\}$ be a basis of $\fn_+^*$ dual to the root vector basis of
$\fn_+$. In this basis, the submanifold $X_w$ is defined by linear equations:
\[
X_w=\{x_\alpha=0,\,\alpha\in w(\Delta_+)\setminus\Delta_+\}.
\]
By definition, the space $\Gamma(U^w, i_{w,+}\cO_{X_w})$ is a module over the Weyl algebra (one
generated by $x_\alpha,\partial/\partial x_\alpha$, $\alpha\in w(\Delta_+)$) with one
generator, ${\bf 1}_w$, and relations
\[
x_\alpha{\bf 1}_w=\partial/\partial x_\beta{\bf 1}_w=0\text{ if }\alpha\in
w(\Delta_+)\setminus\Delta_+,\beta\in w(\Delta_+)\cap\Delta_+.
\]

Similarly, $\Gamma(U^w, \cZ hu_{0}(i_{w,+}\cO_{X_w}))$, is a $\cD^{ch}(U^w)\otimes H_X$-module
with generator ${\bf 1}_w$ and relations
\begin{align}
&\label{descr-wak-w-tw-1} x_{\alpha,n+1}{\bf 1}_w=\partial_{\alpha,n+1}{\bf 1}_w=0;
\text{ if }\alpha\in w(\Delta_+), n\geq 0, \\
\label{descr-wak-w-tw-2} &x_{\alpha,0}{\bf 1}_w=\partial/\partial x_{\beta,0}{\bf 1}_w=0\text{
if }\alpha\in
w(\Delta_+)\setminus\Delta_+,\beta\in w(\Delta_+)\cap\Delta_+,\\
\label{descr-wak-w-tw-char-1}
&(H_X)_{(m)}{\bf 1}_w=0\text{ if }m\in\ZZ
\end{align}

All of this is easy to twist by a character $\nu(z)=\nu_0/z+\nu_{-1}+\cdots$. There still is a
functor, \cite{BB2} sect.2.5.5,
\[
i_{w,+}:\cD^{i^*_w(\cL_{\nu_0})}_{X_w}-{\rm Mod}\rightarrow \cD^{\nu_0}_X-{\rm Mod}.
\]
This gives us a collection of $\cD^{\nu_0}_X$-modules, $i_{w,+}i^*_w\cL_{\nu_0}$, $w\in W$, if
$\nu_0$ is integral. Note that since $X_w$ is affine
\begin{equation}
\label{iso-as - o-mod} i_{w,+}i^*_w\cL_{\nu_0}\stackrel{\sim}{\rightarrow}i_{w,+}\cO_{X_w}
\text{ as }\cO_X-{\rm modules},
\end{equation}
but the actions of $\fg$, on the former induced by (\ref{fg-map-tdo}), on the latter by
(\ref{li-alg-act-flag}), are different.

$\Gamma(U^w, \cZ hu_{\nu(z)}(i_{w,+}\cO_{X_w}))$ is a $\cD^{ch}(U^w)\otimes H_X$-module with
generator ${\bf 1}_{w\circ\nu_0}$ and relations (\ref{descr-wak-w-tw-1}),
(\ref{descr-wak-w-tw-2}), and the following replacement of (\ref{descr-wak-w-tw-char-1})
\begin{equation}
\label{descr-wak-w-tw-char-2} p_{(n)}{\bf
1}_{w\circ\lambda}=p(\nu(z))_{(n)}{\bf 1}_{w\circ\lambda},\; p\in
H_X,
\end{equation}
where $p(\nu(z))_{(n)}$ stands for ${\rm Res}_{z=0}z^{n}p(\nu(z))$.

 By pull-back, each $\Gamma(U^w, \cZ hu_{\nu(z)}(i_{w,+}\cO_{X_w}))$ is a $V_{-h^{\vee}}(\fg)$-
 hence $\ghat$-module. It is quite clear that $\Gamma(U^{{\rm id}}, \cZ hu_{0}(i_{*}\cO_{X_{{\rm id}}}))$
 is precisely the Wakimoto module of critical level and zero central character, $\WW_{\nu(z)=0}$,
 that was introduced in sect.~\ref{def-of-wak-of-crit-1}. We can now see how by passing to the
 twisted CDO we have gained considerable flexibility: denote by
 \begin{equation}
 \label{def-wak-crit-lev-arb-centr-charge}
 \WW_{\nu(z)}=\Gamma(U^{{\rm id}}, \cZ hu_{\nu(z)}(i_{*}\cO_{X_{{\rm id}}}));
 \end{equation}
 these are Wakimoto modules of critical level and central character $\nu(z)$, by construction,
 \cite{FF1,F3}.

 Furthermore, as a quick scan of \cite{F3}, sect.~9.5.1 shows,
 \begin{equation}
 \label{def-wak-crit-lev-arb-centr-charge-tw-w}
 \WW_{\nu(z)}^{w}\stackrel{{\rm def}}{=}\Gamma(U^{w}, \cZ
 hu_{\nu(z)}(i_{w,+}i_w^*\cL_{\nu_0})), \; w\in W,
 \end{equation}
are the so-called $w$-{\it twisted} Wakimoto modules of critical
level and central character $\nu(z)$. (Indeed, used in \cite{F3},
sect.~9.5.1 is the ``$\beta{\rm-}\gamma$-system'' $M_\fg^w$, which
is our $\cD^{ch}(U^{{\rm id}})$ except that the choice of vacuum is
different. To give $\cD^{ch}(U^{{\rm id}})$ a $V_{-h^\vee}(\fg)$-module
structure, the standard morphism $\pi_{{\rm
id}}:V_{-h^\vee}(\fg)\rightarrow\cD^{ch}(U^{{\rm id}})$ is twisted in
\cite{F3}, sect.~9.5.1 by the Tits lifting of $w\in
W\subset\text{Aut}(\fh^*)$ to $\tilde{w}\in\text{Aut}(\fg)$. If
we identify $\cD^{ch}(U^{{\rm id}})$ with $\cD^{ch}(U^{{\rm w}})$
via the same $\tilde{w}\in\text{Aut}(\fg)$, see
sect.~\ref{Feigin-Frenkel-Wakimoto bozonization}, then, by
definition, $\pi_w=\pi\circ \tilde{w}^{-1}$, where $\pi_w$ comes
from (\ref{aff-vetr-alg-cdo-cell-twist}). According to
Theorem~\ref{exist-tw-ch-do-flag} (4), it is these $\{\pi_w\}$ that
`conspire' to define a $V_{-h^\vee}(\fg)$-module structure on
$\cD^{ch,tw}_X$. Now the obvious observation that our choice of
vacuum
(\ref{descr-wak-w-tw-1},\ref{descr-wak-w-tw-2},\ref{descr-wak-w-tw-char-1})
is consistent with the one made in \cite{F3}, sect.~9.5.1 concludes
this little bit of translation.)

\begin{lem}
\label{vanish-higher-delta-fnct}
\begin{align}
\nonumber & H^0(X,\cZ
 hu_{\nu(z)}(i_{w,+}i_w^*\cL_{\nu_0}))=\Gamma(U^{{w}}, \cZ
 hu_{\nu(z)}(i_{w,+}i_w^*\cL_{\nu_0})),\\
 \nonumber
&H^i(X,\cZ
 hu_{\nu(z)}(i_{w,+}i_w^*\cL_{\nu_0}))=0\text{   if  }i>0.
 \end{align}
 \end{lem}

{\em Proof.} As (\ref{sheaf-giving-gr-obj-tw}) shows, $\cZ
 hu_{\nu(z)}(i_{w,+}i_w^*\cL_{\nu_0})$ carries a filtration such that the corresponding
 graded object,
 ${\rm gr}\cZ
 hu_{\nu(z)}(i_{w,+}i_w^*\cL_{\nu_0})$, is an $\cO_X$-module. Furthermore, this graded object is actually a
 push-forward, $i_{w,\ast}\cE$, of a locally free sheaf of $\cO_{X_w}$-modules: to obtain
 this sheaf simply replace in (\ref{sheaf-giving-gr-obj-tw}) $\cT_{X}$ with  $i_w^*\cT_{X}$,
 $\Omega_{X}$ with  $i_w^*\Omega_{X}$, $\cA$ with $S^{\bullet}\cN_{X_w}$, where $\cN_{X_w}$
 is the normal bundle to $X_w$.

 The fact that $X_w$ is affine implies that
 \[
 H^0(X,i_{w,\ast}\cE)=\Gamma(U^{{\rm w}}, i_{w,\ast}\cE),\\
 \nonumber
H^i(X,i_{w,\ast}\cE)=0\text{ if }i>0.
 \]
An application of the standard spectral sequence associated with this filtration gives the
assertion of Lemma~\ref{vanish-higher-delta-fnct} at once. $\qed$

\section{Drinfeld-Sokolov Reduction At The Critical Level}
\label{Drinfeld-Sokolov Reduction At The Critical Level}
\subsection{Categories of  %\boldmath
{$\ghat$}-modules } \label{A category of ghat-modules}
 A triangular decomposition
$\fg=\fn_+\oplus\fh\oplus\fn_-$ determines a triangular
decomposition $\ghat=\hat{\fn}_+\oplus\hat{\fh}\oplus\hat{\fn}_-$,
cf. sect.~\ref{Affine vertex algebras.}, where
$\hat{\fh}=\fh\oplus\CC K$, $\hat{\fn}_{\pm}$ is the preimage of
$\fn_{\pm}$ w.r.t. the evaluation map $\fg[t^{\pm1}]\rightarrow\fg$
defined by letting $t\rightarrow 0$ or $\infty$ resp.

\subsubsection{Definitions of categories }\label{def-of-o-cat}
{\em Define} $\hat{\cO}_k$ to be the category
consisting of $\ghat$-modules $M$ that satisfy the following
conditions:

(1) weight space decomposition: if we let $\hat{\fh}^*_k$ be the
subspace of $\hat{\fh}^*$ defined by the equation $K=k$, then
\begin{equation}
\label{def-weight-space-decomp} M=\oplus_{\mu\in
\hat{\fh}^*_k}M_{\mu},\; M_{\mu}=\{m\in M:\;
hm=\mu(h)m,h\in\hat{\fh}\};
\end{equation}

(2) local finiteness: for each $m\in M$
\begin{equation}
\label{def-loc-finit} {\rm dim}U(\hat{\fn}_+)m<\infty;
\end{equation}

(3) conformal filtration: there exist a family of subspaces
\[\cdots\subset{\rm F}_nM\subset {\rm F}_{n+1}M \cdots\subset M,\;
M=\cup_{n\in\ZZ}{\rm F}_nM,
\]
compatible with the grading of $\ghat$ in that
\begin{equation}
\label{def-filtr-1}\fg\otimes t^m {\rm F}_nM\subset {\rm
F}_{n-m}M,\;{\rm F}_nM=\{0\}\text{ if }n\ll 0
\end{equation}
and locally finitely generated: for each $n\in\ZZ$ there are $m_1,...,m_s$ such that
\begin{equation}
\label{def-loc-fini-gener} {\rm F}_nM\subset\sum_{j=1}^s
U(\ghat)m_j.
\end{equation}

\bigskip
If the {\em level} $k=-h^{\vee}$, the case we shall be interested in
almost exclusively, then we shall use the notation
$\hat{\cO}_{crit}=\hat{\cO}_{-h^{\vee}}$.

Condition (\ref{def-filtr-1}) implies that an object of $\hat{\cO}_k$ is automatically a
$V_k(\fg)$-module. Furthermore, if $k=-h^{\vee}$, then by pull-back, an object of $\hat{\cO}_k$
is automatically a $\fz(\fg)$-module, see Theorem~\ref{exist=of-center-aff-crit}.
One-dimensional irreducible $\fz(\fg)$-modules are nothing but characters, $\chi(z)$, Laurent
series with values in the space dual to the linear span of the generating set
$\{p_1,...,p_r\}$. Given such $\chi(z)=\sum_n\chi_{(n)}z^{-n-1}$, we have a $\fz(\fg)$-module
$\CC_{\chi(z)}$ to be $\CC$ as a vector space with action
$p_i\mapsto\sum_{n}\langle\chi_{(n)},p_i\rangle z^{-n-1}$.

Define $\hat{\cO}_{\chi(z)}$ to be a full subcategory of
$\hat{\cO}_{crit}$ consisting of modules $M$  that satisfy: for each
$m\in M$ and $i$, $1\leq i\leq {\rm rk}\fg$, there is $N$ such that
\begin{equation}
\label{def-of-cat-fix-char} \left(p_{i,
(n)}-\langle\chi_{(n)},p_i\rangle\right)^{N}m=0.
\end{equation}

\subsubsection{Examples of \boldmath{$\ghat$}-modules}
\label{exam-ghat-modul} A rich supply of objects of $\hat{\cO}_{k}$
is obtained by induction: for each $\fg$-module $M$, finite
dimensional or an object of appropriately defined $\cO$-category of
$\fg$-modules, define
\[
{\rm Ind}_{\fg[t]\oplus\CC K}^{\ghat}M,
\]
where $\fg[t]$ operates on $M$ via the evaluation $t\rightarrow 0$
and $K\mapsto k$. In this way we obtain

the Weyl module $\VV_{\lambda,k}={\rm Ind}_{\fg[t]\oplus\CC
K}^{\ghat} V_\lambda$, where $V_\lambda$ is the finite dimensional
simple $\fg$-module with highest weight $\lambda$; note that
$V_{k}(\fg)=\VV_{0,k}$;

the Verma module $\MM_{\lambda,k}={\rm Ind}_{\fg[t]\oplus\CC
K}^{\ghat} M_\lambda$, where $M_\lambda={\rm
Ind}_{\fn_+\oplus\fh}^{\fg}\CC_{\lambda}$, the Verma module over
$\fg$.

The Wakimoto module $\WW_{0,-h^{\vee}}$, see sect.~\ref{Wakimoto
modules.}, belongs to $\hat{\cO}_{crit}$; it is  obtained not so
much by induction as by {\it semi-infinite induction},
\cite{V2,FBZ,F3}.

If $k=-h^{\vee}$, we can introduce {\em restricted} versions, those
obtained by quotienting out by a central character. For
any
$M\in\hat{\cO}_{crit}$ and a
central character $\chi(z)$, define
\[
M_{\chi(z)}=M/{\rm span}\{\left(p_{i,
(n)}-\langle\chi_{(n)},p_i\rangle\right)M\}.
\]
Thus we obtain $\VV_{\chi(z)}$ and $\MM_{\chi(z)}$; these are objects of $\hat{\cO}_{\chi(z)}$.

Twisted Wakimoto modules of critical level, $\WW^w_{\nu(z)}$, which were obtained via
localization in sect.~\ref{examples-of-d-ch-tw-mod},
 are objects of
$\hat{\cO}_{\nu(z)\circ\pi_w}$, where $\nu(z)\circ\pi_w$ stands for the composition of the
vertex algebra morphism $\pi_w|_{\fz(\fg)}$, see (\ref{aff-vetr-alg-cdo-cell-twist}) and the
character $\nu(z):H_X\rightarrow Fields(\CC)$. Note that $\WW^w_{\nu(z)}$ is not necessarily a
quotient of some bigger module by a central character.

\subsubsection{Modules with a Verma flag}
\label{Modules with a Verma flag} We shall say that
$M\in\hat{\cO}_k$ is filtered by Verma modules if it carries a
filtration $G_0M\subset G_1M\subset\cdots$, $\cup_iG_iM=M$, such
that, for each $i$, $G_iM/G_{i-1}M$ is a direct sum of Verma modules
$\MM_{\lambda,k}$, $\lambda\in\fh^*$.

Let $\{v_1,v_2,...\}$ generate $M$ and denote by $V$ the $\hat{\fn_+}\oplus\hat{\fh}$-submodule
of $M$ generated by $\{v_1,v_2,...\}$. Define $P^0\stackrel{{\rm
def}}{=}U(\ghat)\otimes_{\hat{\fn_+}\oplus\hat{\fh}}V$. It is clear that $P^0$ has a filtration
by Verma modules and projects onto $M$:
\[
P^0\rightarrow M.
\]

Continuing in the same vein we obtain, for each $M\in \hat{\cO}_k$ a
resolution
\begin{equation}
\label{resol-by-mod-verma-filtr} \cdots\rightarrow
P^{-j}\rightarrow\cdots\rightarrow P^{-1}\rightarrow
P^{0}\rightarrow M
\end{equation}
by modules with a Verma flag.

Now the locally finite generation condition
(\ref{def-loc-fini-gener}) implies
\begin{lem}
\label{exist-of-res-verm-filt-exhaustive} For each $M\in
\hat{\cO}_k$ and any conformal filtration $\{{\rm F}_nM\}$ there exist
a resolution (\ref{resol-by-mod-verma-filtr}) and conformal
filtrations $\{{\rm F}_nP^{-j}\}$ of all the terms so that

(1) the differential is a morphism of filtered modules;

(2) for each $n$, there is $N$ such that ${\rm F}_mP^{-j}=\{0\}$ for
all $j>N$, $m<n$.
\end{lem}
This is all standard, cf. \cite{RCW}, sect.4, and we omit the details.

\subsubsection{Action of the center and a decomposition into blocks}
\label{Action of the center and  decomposition into blocks} For
generators $p_1,...,p_r$ of $\fz(\fg)$ chosen as in
Theorem~\ref{exist=of-center-aff-crit}, `Fourier  coefficients'
$p_{(n)}$ can be regarded as elements of the completed universal
enveloping algebra $\tilde{U}(\ghat)$, cf. \cite{FBZ}, sect.~4.3.

The conformal weight zero subalgebra $\tilde{U}(\ghat)_0$ has an
ideal $(\tilde{U}(\ghat)U(\fg\otimes\CC[t]t))_0$. We have an obvious
isomorphism
\[
\tilde{U}(\ghat)_0/(\tilde{U}(\ghat)U(\fg\otimes\CC[t]t))_0\stackrel{\sim}{\rightarrow}
U(\fg).
\]

Note that $\CC[p_1,...,p_r]$ is naturally a commutative associative
subalgebra of $\fz(\fg)$ w.r.t. to multiplication $_{(-1)}$.
\begin{lem}
\label{aff-center-fin-center} (\cite{A}) The composite map
\[
\CC[p_1,...,p_r]\rightarrow \tilde{U}(\ghat)_0\rightarrow U(\fg),
p\mapsto p_{0}\, {\rm mod
}\,(\tilde{U}(\ghat)U(\fg\otimes\CC[t]t))_0
\]
has the center $Z(\fg)$ as its image and delivers a commutative
associative algebra isomorphism
$\CC[p_1,...,p_r]\stackrel{\sim}{\rightarrow} Z(\fg)$.
\end{lem}
Note that $p_0$  stands for the {\em conformal weight zero} Fourier
coefficient of the field $p(z)$. For example, since $p_i$ has
conformal weight $d_i+1$, we have
$p_i(z)=\sum_np_{i,n}z^{-n-d_i-1}$, which shows that $p_{i,0}$ is
the coefficient of $z^{-d_i-1}$.

The significance of Lemma~\ref{aff-center-fin-center} is that it
describes the action of the center on Verma modules. To see this,
compose the isomorphism
$\CC[p_1,...,p_r]\stackrel{\sim}{\rightarrow} Z(\fg)$ with the
classical Harish-Chandra isomorphism $\theta:
Z(\fg)\rightarrow\CC[\fh^*]^{W}$ to obtain an isomorphism
\begin{equation}
\label{affine-h-ch-version} \theta_{aff}:
\CC[p_1,...,p_r]\stackrel{\sim}{\rightarrow}\CC[\fh^*]^{W}.
\end{equation}
It follows from the definition of the Verma module
$\MM_{\lambda,-h^{\vee}}$ that for each $f\in\CC[p_1,...,p_r]$,
$f_{n}$ acts on $\MM_{\lambda,-h^{\vee}}$ as 0 if $n>0$ and as
multiplication by $\langle\theta_{aff}(f),\lambda\rangle$ if $n=0$.

The same applies to any highest weight module, and, since any object
of $\hat{\cO}_{crit}$ has a filtration by highest weight modules
(alternatively, use Lemma~\ref{exist-of-res-verm-filt-exhaustive}),
one obtains a block decomposition
\begin{equation}
\label{block-decomp}
\hat{\cO}_{crit}=\oplus_{[\lambda]\in\fh^*/W}\hat{\cO}_{crit}^{[\lambda]},
\end{equation}
where $[\lambda]=W\circ\lambda$ and  $\hat{\cO}_{crit}^{[\lambda]}$
is defined to be the full subcategory of modules such that
$(f_{0}-\theta_{aff}(f)(\lambda))$ acts locally nilpotently  for
each $f\in\CC[p_1,...,p_r]$.

Introduce a polynomial algebra
\begin{equation}
\label{def-of-z-min} \cZ_-=\CC[p_{i,n},1\leq i\leq {\rm rk}\fg,n<0].
\end{equation}
Its is clear that $f_{n}$, $f\in\fz(\fg)$, $n>0$, acts locally
nilpotently on each $M\in \hat{\cO}_{crit}$. The action of $f_{0}$
being described by (\ref{block-decomp}), what one needs to describe
the action of the entire $\fz(\fg)$ is the action of $\cZ_-$. Here
is an example.
\begin{thm} (\cite{F3}, Theorem~9.5.3)
\label{centre-on-verma-frenk} $\MM_{\lambda, -h^{\vee}}$ is a free
$\cZ_-$-module. Furthermore, $\text{End}_{\widehat{\fg}}(\MM_{\lambda, -h^{\vee}})\stackrel{\sim}{\rightarrow}\cZ_-$.
\end{thm}

\subsubsection{Gradings and character formulas}
\label{Gradings and character formulas} A conformal filtration that $M\in\hat{\cO}_k$ carries
by definition, see sect.~\ref{def-of-o-cat}, is not unique. But if $M$ is a Verma module or its
quotient, then there are obvious choices: pick an $m$, declare ${\rm F}_nM=\{0\}$ if $n<m$, let
${\rm F}_mM$ contain $v$,  a highest weight vector, and define ${\rm
F}_{n+m}M=U(\hat{\fn}_-)_{\geq -n}v$; here $U(\hat{\fn}_-)_{\geq -n}$ stands for the subspace of
$U(\hat{\fn}_-)$ spanned by $x_1\otimes t^{-n_1}\cdots x_l\otimes t^{-n_l}$ with all $n_j\geq
0$ and $\sum_jn_j\leq n$.

 If $M$ is a
highest weight module with filtration, we shall tacitly assume that ${\rm F}_nM=\{0\}$ if $n<0$, and that the highest weight belongs to ${\rm F}_0M$. We shall denote by $M[m]$ this same $M$ with
filtration  shifted by $m$, as in the paragraph above. This explains
the assertion (2) of Theorem~\ref{main-theorem-intro}.

In any case, given a filtration ${\rm F}_\bullet M$ of $M$, the
graded object, $\operatorname{gr}_{{\rm F}}M$, is also a $\ghat$-module, the
action being defined to be that on symbols. Furthermore, the
canonical grading of $\operatorname{gr}_{{\rm F}}M$ makes it into a graded
$\ghat$-module.

This can be made a little more explicit: extend $\ghat$ to
$\ghat_{ext}$ by adjoining, as usual, a {\it degree derivation} $D$
so that
\begin{align*}
[D,x\otimes t^n]=-nx\otimes t^n.
\end{align*}
 Letting $D$ act on ${\rm
F}_nM/{\rm F}_{n-1}M$ as multiplication by $n$, we make $\operatorname{gr}M$ into a
$\ghat_{ext}$-module. More generally, we shall call $M\in\hat{\cO}_k$ {\em graded} if
the action of $\ghat$ extends to that of $\ghat_{ext}$ so that the action of $D$ is diagonalizable. We have the weight space decomposition
\[
M=\bigoplus_{n\in(\CC D)^*}M_n.
\]

This is reflected in the following $q$-{\em
dimension} formula, written down in case of an arbitrary filtered module,
\begin{equation}
\label{q-dim-defin} \dim_qM=\sum_{n\gg -\infty}q^n\dim\operatorname{gr}M_n.
\end{equation}
Note that it makes sense only if $\dim\operatorname{gr}M_n<\infty$. For
example,
\begin{equation}
\label{q-dim-of-weyl} {\rm dim}_q\VV_{\lambda,k}={\rm
dim}V_{\lambda}\prod_{j=1}^{\infty}(1-q^j)^{-{\rm dim}\fg}.
\end{equation}

Invoking the semi-simplicity of the action of $\fh$, we now obtain a weight space decomposition of an arbitrary graded module:
\[
M=\bigoplus_{\alpha\in(\hat{\fh}\oplus\CC D)^*}M_{\alpha},
\]
and refine (\ref{q-dim-defin}) by  defining the {\it
formal character} as usual
\begin{equation}
\label{what - the-char-form-is} {\rm ch} M\stackrel{{\rm def}}{=}{\rm }{\rm ch}( {\rm
gr}M)\stackrel{{\rm def}}{=}\sum_{\alpha\in(\hat{\fh}\oplus\CC D)^*}e^{\alpha}\dim\operatorname{gr}M_{\alpha}.
\end{equation}
Note that, as it follows from the definition of $\hat{\cO}_k$, sect.~\ref{def-of-o-cat}, $\dim\operatorname{gr}M_{\alpha}<\infty$.

For example, the Verma module $\MM_{\lambda,k}$ is actually graded and
\begin{equation}
\label{char-of-verma-mod} {\rm
ch}\MM_{\lambda,k}=e^{\lambda}\prod_{\alpha\in\hat{\Delta}_+}(1-e^{-\alpha})^{-1},
\end{equation}
where $\hat{\Delta}_+$ is the set of positive roots of $\widehat{\fg}$.

The freeness result, Theorem~\ref{centre-on-verma-frenk},    implies
\begin{equation}
\label{char-of-verma-mod-restr} {\rm
ch}\MM_{\mu(z)}=e^{\lambda}\prod_{\alpha\in\Delta_+}(1-e^{-\alpha})^{-1}\times
\prod_{n=1}^{+\infty}\prod_{\alpha\in\Delta_+}(1-e^{-\alpha-n\delta})^{-1}\times
\prod_{\alpha\in\Delta_-}(1-e^{-\alpha-n\delta})^{-1},
\end{equation}
where $\mu(z)$ and $\lambda$ must be compatible, of course: $\mu(z)=\lambda/z+\mu_1+\mu_2z+\cdots.$

Some modules are graded by the very construction, for example
 $w$-twisted Wakimoto modules with {\it homogeneous} character $\nu(z)=\nu_0/z$,
and so the formal character of $\WW^w_{\nu_0/z}$ does not require
specification of a filtration. On the other hand, it is quite clear
how to define a filtration on an arbitrary $w$-twisted Wakimoto
module so that its formal character coincide with that of
$\WW^w_{\nu_0/z}$. The result is, cf. \cite{F3}, formula (9.5.4),
\begin{equation}
\label{char-of-wak-mod-restr} {\rm
ch}\WW^w_{\nu(z)}=e^{w\circ\nu_0}\prod_{\alpha\in\Delta_+}(1-e^{-\alpha})^{-1}\times
\prod_{n=1}^{+\infty}\prod_{\alpha\in\Delta_+}(1-e^{-\alpha-n\delta})^{-1}\times
\prod_{\alpha\in\Delta_-}(1-e^{-\alpha-n\delta})^{-1}.
\end{equation}
This, of course, coincides with the character of the restrticted
Verma module (\ref{char-of-verma-mod-restr}) with highest weight
$\lambda=w\circ\nu_0$.

\subsubsection{ }
\label{duality-for-graded-modules}
Since $M\in \hat{\cO}_k$  being graded implies  $\dim\operatorname{gr}M_{\alpha}<\infty$,  we define its dual
\[
M^c=\bigoplus_{\alpha\in(\hat{\fh}\oplus\CC D)^*}M_{\alpha}^*,
\]
with the action of $\ghat_{ext}$ determined by
$\langle x\phi,m\rangle=\langle \phi,\omega(x)m\rangle$, $x\in\ghat_{ext}$, $\phi\in M^c$,
$m\in M$, and $\omega:\ghat_{ext}\rightarrow\ghat_{ext}$ is the canonical antiinvolution that sends $\fg_\alpha\otimes t^n$ to $\fg_{-\alpha}\otimes t^{-n}$, $\alpha\in\Delta$
being a root of $\fg$.

It is obvious that the assignment $M\mapsto M^c$ is a contragredient exact functor
on the full subcategory of $\hat{\cO}_k$ consisting of graded modules.

\subsection{The Drinfeld-Sokolov reduction}
\label{The Drinfeld-Sokolov reduction}

\subsubsection{Definition}
\label{The Drinfeld-Sokolov reduction Definition and main results}

Consider the vertex superalgebra $V(\fn_+)\otimes Cl(\fn_+)$, cf. sect.~\ref{A super-version:
Clifford algebra.} and \ref{Affine vertex algebras.}. Let $\{e^\alpha,\alpha\in\Delta_+\}$ be a
root vector basis of $\fn_+$, $\{\phi_\alpha,\alpha\in\Delta_+\}$ its copy, albeit with changed
parity, that appears inside $Cl(\fn_+)$, $\{\phi_\alpha^*,\alpha\in\Delta_+\}\subset\fn_+^*$
the basis dual to the latter.

Let $\{c^{\alpha\beta}_\gamma\}$ be the structure constants so that
\[
[e^\alpha,e^\beta]=\sum_\gamma c^{\alpha\beta}_\gamma e^\gamma.
\]
The following  elements of $V(\fn_+)\otimes Cl(\fn_+)$ are of
importance
\[
Q_{st}=\sum_{\alpha}e^\alpha\otimes\phi_\alpha^*-
\frac{1}{2}\sum_{\alpha,\beta,\gamma}c^{\alpha,\beta}_\gamma\phi_{\alpha(-1)}^*(\phi_{\beta(-1)}^*\phi_{\gamma}),\
\chi=\sum_{i=1}^{r}\phi_{\alpha_i}^*, Q_{DS}=Q_{st}+\chi,
\]
where $\{\alpha_1,...,\alpha_r\}$ is the set of simple roots. Define
\[
d_{st}=(Q_{st})_{(0)},\ d=d_{st}+\chi_{(0)}.
\]
Since $d_{st}^2=\chi_{(0)}^2=[d_{st},\chi_{0)}]=0$, there arise 3
differential graded vertex algebras, $(V(\fn_+)\otimes
Cl(\fn_+),d_{st})$, $(V(\fn_+)\otimes Cl(\fn_+),\chi_{(0)})$,
$(V(\fn_+)\otimes Cl(\fn_+),d)$, with grading defined by setting
${\rm deg}\phi^*_\alpha=1$, ${\rm deg}\phi_\alpha=-1$.

Furthermore, if $M$ is a $V(\fn_+)$-module, then $(M\otimes
Cl(\fn_+),d_{st})$, $(M\otimes Cl(\fn_+),\chi_{(0)})$, $(M\otimes
Cl(\fn_+),d)$ are differential graded modules over their respective
differential graded vertex algebras. To emphasize the fact that all
these are to be treated as complexes, we shall change the notation
and write $C^{\infty/2+\bullet}(L\fn_+, M)$ instead of $M\otimes
Cl(\fn_+)$. The corresponding cohomology will be denoted as follows:
$H^{\infty/2+\bullet}(L\fn_+, M)$,
$H_{\chi_{(0)}}^{\infty/2+\bullet}(L\fn_+, M)$,
$H^{\infty/2+\bullet}_{DS}(L\fn_+, M)$.
\begin{rem}
\label{rem-on-dr-sok-funct}
$H^{\infty/2+\bullet}(L\fn_+, M)$,
and
$H^{\infty/2+\bullet}_{DS}(L\fn_+, M)$ can be related  as follows.
One can regard $\chi$ as a (Drinfeld-Sokolov)character, $\chi:L\fn_+\rightarrow\CC$, which sends $e^{\alpha_i}\otimes t^{-1}$ to
1 and the rest of $\{e^{\alpha}\otimes t^{n}\}$ to 0. Denote by $\CC_\chi$ the corresponding 1-dimensional $L\fn_+$-module. Then by definition
$H^{\infty/2+\bullet}_{DS}(L\fn_+, M)=H^{\infty/2+\bullet}(L\fn_+, M\otimes\CC_\chi)$.
\end{rem}
\bigskip

If, in addition, $M\in\hat{\cO}_{crit}$ and is regarded as an $L\fn_+$-module via pull-back,
then each of the three series of cohomology groups above is a $\fz(\fg)$-module, because the
action of the center commutes with that of $L\fn_+$. Thus we obtain three series of functors
\begin{equation}
\label{def-of-dr-sok-red-form} H^{\infty/2+\bullet}(L\fn_+, ?),
H_{\chi_{(0)}}^{\infty/2+\bullet}(L\fn_+, ?),
H^{\infty/2+\bullet}_{DS}(L\fn_+, ?):\; \hat{\cO}_{crit}\rightarrow
\fz(\fg){\rm -Mod}.
\end{equation}
Each of these functors makes sense away from the critical level.
For example, if we let $\cW_k\stackrel{{\rm
def}}{=}H^{\infty/2+\bullet}_{DS}(L\fn_+, V_k(\fg))$, then we obtain
\[
H^{\infty/2+\bullet}_{DS}(L\fn_+, ?):\; \hat{\cO}_{crit}\rightarrow
\cW_k{\rm -Mod}.
\]

 All of this is well known, of course:
$H^{\infty/2+\bullet}(L\fn_+, ?)$ was introduced by Feigin in
\cite{Feig}, the entrance point of the BRST business in mathematics,
$H^{\infty/2+\bullet}_{DS}(L\fn_+, ?)$ is the Drinfeld-Sokolov
reduction functor, proposed by Feigin and Frenkel as a tool to
define $\cW_k$, the celebrated $W$-algebra, see \cite{FBZ,F3} and
references therein. The Drinfeld-Sokolov reduction functor has been
studied in \cite{A2,A4,A} in a more general setting. A thorough
analysis of the functors (\ref{def-of-dr-sok-red-form}) has been
carried out recently by Frenkel and Gaitsgory \cite{FG1,FG2,FG3}.

\subsubsection{Torus action, grading, twisted grading and filtration}
\label{Twisted grading and filtration dr-sok} Note that $Cl(\fn_+)$
carries a natural grading determined by the condition that the
degree of  $\phi_{\alpha,n}$ and $\phi_{\alpha,n}^*$ be $-n$. If
$M\in\hat{\cO}_{crit}$ is graded, i.e., carries an action of the
extended $\ghat_{ext}$ with diagonalizable $D$, cf.
sect.~\ref{Gradings and character formulas}, then the entire complex
$C^{\infty/2+\bullet}(L\fn_+,M)$ acquires a grading, that of the
tensor products of graded spaces, so that $d_{st}$ has degree 0.
Therefore, this grading descends on $H^{\infty/2+\bullet}(L\fn_+,
M)$.

This can be refined by noting that we can extend an action not only of $D$ but of the entire
$L\fh$ to $C^{\infty/2+\bullet}(L\fn_+,M)$. Consider the linear map
\begin{equation}
\label{def-of-cartan-action-1} \fh\rightarrow
C^{\infty/2+\bullet}(L\fn_+,V_{-h^{\vee}}(\fh)),\; h\mapsto h\otimes{\bf
1}-\sum_{\alpha\in\Delta_+}\alpha(h){\bf 1}\otimes\phi^*_{\alpha,(-1)}\phi_{\alpha}.
\end{equation}
It is easy to see that it defines a vertex algebra morphism, cf.the end of sect.~\ref{Affine
vertex algebras.},
\begin{equation}
\label{def-of-cartan-action-2} V_0(\fh)\rightarrow
C^{\infty/2+0}(L\fn_+,V_{-h^{\vee}}(\fg)),
\end{equation}
where the central charge has got shifted: $-h^{\vee}(.,.)$ has been replaced with $0$.

 Therefore,
$C^{\infty/2+\bullet}(L\fn_+,V_{-h^{\vee}}(\fg))$, $M\in \hat{\cO}_{crit}$, is a
$V_0(\fh)$-module, and it is easy to see that $d_{st}$ is a $V_0(\fh)$-module morphism. In
particular, if $M$ is graded and $M=\oplus_{\alpha\in\hat{\fh}^*}M_\alpha$, then
$H^{\infty/2+\bullet}(L\fn_+,M)$ is also graded:
\begin{equation}
\label{Pweight-space-decomp-of-semi-inf}
 H^{\infty/2+\bullet}(L\fn_+,M)=\oplus_{\alpha\in\hat{\fh}^*}
 H^{\infty/2+\bullet}(L\fn_+,M)_\alpha.
 \end{equation}

None of this carries over to the Drinfeld-Sokolov case, because $\chi_{(0)}$ does not preserve
either of the gradings introduced.
 To rescue the situation -- partially -- denote by $\tilde{D}\in{\rm
End}C^{\infty/2+\bullet}(L\fn_+,M)$ the operator whose eigenspace
decomposition coincides with the grading just discussed and
introduce $\hat{\rho}^\vee\in{\rm
End}C^{\infty/2+\bullet}(L\fn_+,M)$, the operator whose eigenvalues
are negative those of the half-sum of positive coroots
$\rho^{\vee}$. Specifically, we demand that
\[
[\hat{\rho}^\vee,e_{\alpha,n}]=-\rho^{\vee}(\alpha)e_{\alpha,n},
[\hat{\rho}^\vee,\phi_{\alpha,n}]=-\rho^{\vee}(\alpha)\phi_{\alpha,n},
[\hat{\rho}^\vee,\phi^*_{\alpha,n}]=\rho^{\vee}(\alpha)\phi^*_{\alpha,n}.
\]
Now notice \cite{FBZ,A} that $\chi_{(0)}$, hence
$d=d_{st}+\chi_{(0)}$, commutes with $\tilde{D}+\hat{\rho}^\vee$,
and the eigenvalues  of the latter provide a {\it twisted grading}
of the Drinfeld-Sokolov reduction $H^{\infty/2+\bullet}_{DS}(L\fn_+,
M)$. This allows us to define, cf. (\ref{q-dim-of-weyl}), the notion
of $q$-{\it dimension}
\begin{equation}
\label{def-q-dim-dr-sok} {\rm
dim}_qH^{\infty/2+\bullet}_{DS}(L\fn_+, M)=\sum_{n\in\ZZ}q^n{\rm
dim}H^{\infty/2+\bullet}_{DS}(L\fn_+, M)_n.
\end{equation}

Similarly, if $M$ is filtered, then this filtration extends to $C^{\infty/2+\bullet}(L\fn_+,M)$
and then descends to $H^{\infty/2+\bullet}_{DS}(L\fn_+, M)$ automatically. A spectral sequence
arises
\begin{equation}
\label{spectr-from-stand-to-ds}
\{E^{\bullet\bullet}_r,d_r\}\Rightarrow
H^{\infty/2+\bullet}_{DS}(L\fn_+, M)\text{ s.t.
}E^{\bullet\bullet}_1,d_1=H^{\infty/2+\bullet}(L\fn_+, M).
\end{equation}

\subsubsection{Drinfeld-Sokolov reduction as a derived functor}

Note that each {\it singular vector}, i.e., an element of
$M^{\hat{\fn}_+}$, is a cocycle:
\begin{equation}
\label{sing-vect-cocycle} M^{\hat{\fn}_+}\hookrightarrow
Z^{\infty/2+0}(L\fn_+,M), m\mapsto m\otimes{\bf 1};
\end{equation}
this is true of any differential, $d$ or $d_{st}$, and is an obvious
consequence of the definition.

\begin{thm}
\label{dr-sok-as-derived}

(1) If $P\in\hat{\cO}_{k}$ is filtered by Verma modules,
then  $H_{DS}^{\infty/2+i}(L\fn_+,P)=0$ for all
$i\neq0$ \cite{A2}. Furthermore, if $k=-h^\vee$ and all Verma modules appearing in
the composition series of $P$ have regular integral highest weight, then
$H_{DS}^{\infty/2+0}(L\fn_+,P)$ is a free $\cZ_-$-module (see
(\ref{def-of-z-min})) on generators $[v]={\rm class}\,v$, where $v$
varies over the set of highest weight vectors of Verma modules that
appear in the associated graded of $P$.

(2) If $i>0$, then $H^{\infty/2+i}_{DS}(L\fn_+, M)=0$ for all
$M\in\hat{\cO}_{k}$.

(3) The functor $H^{\infty/2+0}_{DS}(L\fn_+, ?)$ is right exact.

(4) The class of modules carrying a filtration by Verma modules is
adapted to the functor $H^{\infty/2+0}_{DS}(L\fn_+, ?)$.
\end{thm}
Let us discuss some consequences. Items (3) and (4) allow us to
define, as usual \cite{GelMan}, the derived functors $L^i
H^{\infty/2+0}_{DS}(L\fn_+, ?)$.
\begin{cor}
\label{iso-with-der-funct} The functors $H^{\infty/2-i}_{DS}(L\fn_+,
?)$ and $L^i H^{\infty/2+0}_{DS}(L\fn_+, ?)$ are isomorphic.
\end{cor}

\begin{cor}If $\mu_0$ is regular integral, then
\label{dr-sok-of-restr-verma}
\[
H^{\infty/2+i}_{DS}(L\fn_+, \MM_{\mu(z)})=\left\{
\begin{array}{ll}
\CC[v_{\mu_0}]&\text{ if }i=0\\
0&\text{ otherwise, }
\end{array}
\right.
\]
where $[v_{\mu_0}]$ is the cohomology class of a highest weight
 vector $v_{\mu_0}$.
\end{cor}

{\em Proof of Corollary~\ref{iso-with-der-funct} .} Pick a
resolution of $M$
\[
P^{\bullet}:\;\cdots\rightarrow P^{-j}\rightarrow\cdots\rightarrow
P^{-1}\rightarrow P^{0}\rightarrow M
\]
by modules with Verma filtration. By definition, $L^{\bullet}
H^{\infty/2+0}_{DS}(L\fn_+, ?)$ is the cohomology of the complex
$H^{\infty/2+0}_{DS}(L\fn_+, P^\bullet)$. That the latter complex
also computes $H^{\infty/2-\bullet}_{DS}(L\fn_+, M)$ is derived from
Theorem~\ref{dr-sok-as-derived} (1) and (3) by a standard
argument; it is based on the long exact sequences of
various $H^{\infty/2-\bullet}_{DS}(L\fn_+, ?)$ that are associated
to the short exact sequences
\[
0\rightarrow \cZ^0\rightarrow P^0\rightarrow M\rightarrow 0,
\]
\[
0\rightarrow \cZ^{-i-1}\rightarrow P^{-i-1}\rightarrow
\cZ^{-i}\rightarrow 0,\;i\geq 0.
\]
It is that argument which allows to compute the cohomology of a
sheaf via its $\bar{\partial}$-resolution. We leave the details to
the interested reader. $\qed$

{\em Proof of Corollary~\ref{dr-sok-of-restr-verma}.} For a
commutative algebra $A$, a collection of elements
$\underline{a}\subset A$, and an $A$-module $E$, denote by
$K^\bullet(A,\underline{a};E)$ the corresponding Koszul complex.

Let now $\underline{z}=\{p_{in}-p_i(\nu_n),\; n<0,1\leq i\leq {\rm
rk}\fg\}\subset\cZ_-$. This collection being regular,
Theorem~\ref{centre-on-verma-frenk} implies that
$K^\bullet(\cZ_-,\underline{z};\MM_{\nu_{0},-h^{\vee}})$ is a
resolution of $\MM_{\mu(z)}$. Due to
Corollary~\ref{iso-with-der-funct}, the complex
$H^{\infty/2+0}(L\fn_+,
K^\bullet(\cZ_-,\underline{z};\MM_{\nu_{0},-h^{\vee}}))$ computes
$H^{\infty/2-\bullet}(L\fn_+,\MM_{\mu(z)})$.
Theorem~\ref{dr-sok-as-derived}(1) says that the latter complex is
nothing but $K^\bullet(\cZ_-,\underline{z};\cZ_-)$, the Koszul
resolution of $\cZ/\langle\underline{z}\rangle=\CC$, where by
$\langle\underline{z}\rangle$ we have denoted the ideal generated by
$\underline{z}$. $\qed$

\subsubsection{Proof of Theorem~\ref{dr-sok-as-derived} (1)}
\label{proof-of-theor-dr-sok-as-derived 1} Except
for the freeness assertion in the case where $k=-h^\vee$, this item is proved in \cite{A2}. Let us review the details.

First of all,  \cite{A2}, Theorem~5.7 and Remark~5.8,
\begin{equation}
\label{arakawa-coho-verma} {\rm dim}_q
H^i_{DS}(L\fn_+,\MM_{\lambda,-h^{\vee}})=\left\{
\begin{array}{lll}0&\text{ if }&i\neq 0\\
q^{\langle\lambda,D-\rho^{\vee}\rangle }
\prod_{j=1}^{\infty}(1-q^j)^{-{\rm rk}\fg}&\text{ if }& i= 0,
\end{array}\right.
\end{equation}
cf. (\ref{def-q-dim-dr-sok}); note that our
conformal grading convention is different from that in \cite{A2}.
 Next, if $P$ is an arbitrary module with Verma
filtration, then there arises the standard spectral sequence
associated to this filtration, which is easily seen to converge to
$H^{\infty/2+\bullet}_{DS}(L\fn_+,P)$. The vanishing result
(\ref{arakawa-coho-verma}) implies that (a) the spectral sequence
collapses in the first term, (b) $H^{\infty/2+i}_{DS}(L\fn_+,P)=0$ if $i\neq
0$, (c) $H^{\infty/2+0}_{DS}(L\fn_+,P)$ is filtered so that the corresponding
graded object is a direct sum of various
$H^{\infty/2+0}_{DS}(L\fn_+,\MM_{\lambda,-h^{\vee}})$, one for each Verma
module occurring in the Verma composition series of $P$. This proves
\begin{lem}
\label{dr-sok-of-mod-with-verma}
If $P$ is a module with Verma filtration so that $\operatorname{gr}P=\oplus_{\alpha\in A}\MM_{\lambda_\alpha,-h^\vee}$, $A$ being an index set, then
\[
\operatorname{gr}H^{\infty/2+i}_{DS}(L\fn_+, P)=\left\{
\begin{array}{ll}
\oplus_{\alpha\in A}H^{\infty/2+0}_{DS}(L\fn_+, \MM_{\lambda_\alpha,-h^\vee})&\text{ if }i=0\\
0&\text{ otherwise, }
\end{array}
\right.
\]
\end{lem}

As to the freeness assertion, let us start off by making an informal
remark. If $v_{\lambda}$ is a highest weight vector of
$\MM_{\lambda,k}$, then  according to (\ref{sing-vect-cocycle}) it
determines a cohomology class, $[v_\lambda]$. Note that its twisted
degree, sect.~\ref{Twisted grading and filtration dr-sok}, is
precisely $\langle\lambda,D-\rho^{\vee}\rangle$.  Since the
$q$-dimension of $\cZ_-$ is clearly
$\prod_{j=1}^{\infty}(1-q^j)^{-{\rm rk}\fg}$, the second line of
(\ref{arakawa-coho-verma}) is a strong indication that $[v_\lambda]$
freely generates $H^0_{DS}(L\fn_+,\MM_{\lambda,-h^{\vee}})$. Furthermore,
if each $H^{\infty/2+i}_{DS}(L\fn_+, \MM_{\lambda_\alpha,-h^\vee})$ is a free
$\cZ_-$-module on one generator, then the following refinement of Lemma~\ref{dr-sok-of-mod-with-verma}
is valid:
\[
H^{\infty/2+i}_{DS}(L\fn_+, P)=\oplus_{\alpha\in A}\cZ_-,
\]
because a filtration by free modules splits. It remains then to prove
\begin{thm}
\label{freeness-ds-for-verma}
If $\lambda$ is a regular integral weight, then $H^{\infty/2+0}_{DS}(L\fn_+, \MM_{\lambda,-h^\vee})$ is a free $\cZ_-$-module on one generator $[v_\lambda]$.
\end{thm}
Since the rest of Theorem~\ref{dr-sok-as-derived} is independent of this result, we will
postpone proving Theorem~\ref{freeness-ds-for-verma} until sect.\ref{proof-of thm-on-free-ds-verma}.

\subsubsection{Proof of Theorem~\ref{dr-sok-as-derived} (2)}
Given  $M\in\hat{\cO}_{k}$, consider a resolution of $M$ by
modules  with a Verma filtration
(\ref{resol-by-mod-verma-filtr}). Letting $Z^{-j}=\text{Ker}(P^{-j}\rightarrow P^{-j+1})$
we obtain a collection of short exact sequences
\[
0\rightarrow Z^{0}\rightarrow P^{0}\rightarrow M\rightarrow 0,\;
0\rightarrow Z^{-j}\rightarrow P^{-j}\rightarrow Z^{-j+1}\rightarrow 0,\; j>0.
\]
By virtue of the vanishing result in Lemma~\ref{dr-sok-of-mod-with-verma}, an application of the long exact
sequence of cohomology gives, for each $i>0$ and $n\in\ZZ$, a chain
of isomorphisms
\[
\begin{aligned}
&{\rm F}_nH^{\infty/2+i}(L\fn_+,M)\stackrel{\sim}{\rightarrow}
{\rm F}_nH^{\infty/2+i+1}(L\fn_+,Z^0)\stackrel{\sim}{\rightarrow}\\
&{\rm F}_nH^{\infty/2+i+2}(L\fn_+,Z^{-1})
\stackrel{\sim}{\rightarrow}\cdots \stackrel{\sim}{\rightarrow}{\rm
F}_nH^{\infty/2+j+1}(L\fn_+,Z^{-j})\stackrel{\sim}{\rightarrow}\cdots,
\end{aligned}
\]
where ${\rm F}_n$ denotes the n-th conformal filtration component, see sect.\ref{def-of-o-cat}.
\begin{sloppypar}
Making sure that the resolution $P^\bullet$ is one from
Lemma~\ref{exist-of-res-verm-filt-exhaustive}, we see that ${\rm
F}_nH^{\infty/2+j+1}(L\fn_+,Z^{-j})=\{0\}$ if $j\gg 0$. $\qed$
\end{sloppypar}

\subsubsection{Proof of Theorem~\ref{dr-sok-as-derived} (3)}
If $0\rightarrow A\rightarrow B\rightarrow C\rightarrow0$ is exact, we obtain the long exact
sequence of cohomology
\[
\begin{aligned}
\cdots\rightarrow &H^{\infty/2+0}_{DS}(L\fn_+,A)\rightarrow
H^{\infty/2+0}_{DS}(L\fn_+,B)\rightarrow
H^{\infty/2+0}_{DS}(L\fn_+,C)\rightarrow\\
&H^{\infty/2+1}_{DS}(L\fn_+,A)\rightarrow\cdots .
\end{aligned}
\]
The right exactness follows, because by virtue of Theorem~\ref{dr-sok-as-derived} (2),
$H^{\infty/2+1}_{DS}(L\fn_+,A)=0$.

\subsubsection{Proof of Theorem~\ref{dr-sok-as-derived} (4)} Being adapted means \cite{GelMan}
that (a) each module has a resolution by modules with Verma
filtration, (b) if a complex $P^\bullet$ consisting of modules with
Verma filtration is exact, then
$H^{\infty/2+0}_{DS}(L\fn_+,P^\bullet)$ is also exact. Item (a) is
the assertion of Lemma~\ref{exist-of-res-verm-filt-exhaustive}. Item
(b) is a standard consequence of Theorem~\ref{dr-sok-as-derived} (1)
and (3): present an exact
sequence
\[
\cdots \rightarrow P^{-n}\rightarrow\cdots\rightarrow P^0\rightarrow 0
\]
as a chain of short exact sequences
\[0\rightarrow Z^{-1}\rightarrow P^{-1}\rightarrow P^{0}\rightarrow 0,\;
0\rightarrow Z^{-j-1}\rightarrow P^{-j-1}\rightarrow Z^{-j}\rightarrow 0,j>0.
\]
Then an induction on $j$, using Theorem~\ref{dr-sok-as-derived} (1)
and (3),
 will show that
\[
\cdots \rightarrow H^{\infty/2+0}_{DS}(L\fn_+,P^{-n})\rightarrow\cdots\rightarrow
H^{\infty/2+0}_{DS}(L\fn_+,P^0)\rightarrow 0
\]
is the composition of the short exact sequences
\[
0\rightarrow H^{\infty/2+0}_{DS}(L\fn_+,Z^{-1})\rightarrow
H^{\infty/2+0}_{DS}(L\fn_+,P^{-1})\rightarrow H^{\infty/2+0}_{DS}(L\fn_+,P^{0})\rightarrow 0,
\]
\[
0\rightarrow H^{\infty/2+0}_{DS}(L\fn_+,Z^{-j-1})\rightarrow
H^{\infty/2+0}_{DS}(L\fn_+,P^{-j-1})\rightarrow H^{\infty/2+0}_{DS}(L\fn_+,Z^{-j})\rightarrow
0,j>0,
\]
and is, therefore, also exact.
\subsection{Proof of Theorem~\ref{freeness-ds-for-verma}}
\label{proof-of thm-on-free-ds-verma}

\subsubsection{ }
\label{two-reductions} We begin with a series of 2 reductions, to be followed by
2 ^^ ^^ approximations." The first reduction is the content of
\begin{lem}
\label{reduction-one} Let $\lambda$ be a regular integral weight. The following conditions are equivalent:

(i) $H^{\infty/2+0}_{DS}(L\fn_+,\MM_{\lambda, -h^\vee})$ is a rank 1 free $\cZ_-$-module
generated by $[v_\lambda]$;

(ii) $H^{\infty/2+0}_{DS}(L\fn_+,\MM_{\lambda, -h^\vee})$ is a cyclic $\cZ_-$-module
generated by $[v_\lambda]$;

(iii) $H^{\infty/2+0}_{DS}(L\fn_+,\MM_{\lambda, -h^\vee})$ is a rank 1 free $\cZ_-$-module
generated by $[v_\lambda]$ if $\lambda$ is dominant.
\end{lem}
{\em Proof.}

Equivalence $(i)\Leftrightarrow (ii)$. The implication $(i)\Rightarrow (ii)$ being clear, it
suffices to show $(i)\Leftarrow (ii)$. Item $(ii)$ implies a surjection of $\cZ_-$-modules
$\cZ_-\rightarrow H^{\infty/2+0}_{DS}(L\fn_+,\MM_{\lambda, -h^\vee})$. On the other hand, formula (\ref{arakawa-coho-verma}) shows that $\dim_q\cZ_-=\dim_q  H^{\infty/2+0}_{DS}(L\fn_+,\MM_{\lambda, -h^\vee})$, cf. a discussion in sect.\ref{proof-of-theor-dr-sok-as-derived 1}. Hence this surjection is an isomorphism.

Equivalence $(i)\Leftrightarrow (iii)$. We only need the implication $(i)\Leftarrow(iii)$, and it suffices to show that the class $[v_{w\circ\lambda}]$, $w\in W$, is not annihilated by any non-zero element of $\cZ_-$. (Consider the map $\cZ_-\rightarrow H^{\infty/2+0}_{DS}(L\fn_+,\MM_{w\circ\lambda, -h^\vee})$, $P\mapsto P[v_{w\circ\lambda}]$, and use the character formula (\ref{arakawa-coho-verma}), as we have just done.)

Consider a reduced expression $w=s_{i_k}s_{i_{k-1}}\cdots s_{i_1}$. It is known that any
morphism of Verma modules (over $\fg$) $M_{w\circ\lambda}\rightarrow M_\lambda$ factors out in the composition of morphisms $M_{s_{i_{j}}\circ\lambda_{j-1}}\rightarrow M_{\lambda_{j-1}}$ determined by the assignment $v_{s_{i_{j}}\circ\lambda_{j-1}}\mapsto
f_{i_{j}}^{<\lambda,\alpha_{i_j}^\vee>+1}v_{\lambda_{j-1}}$; here $\lambda_{j-1}=(s_{i_{j-1}}\cdots s_{i_1})\circ\lambda$, $f_{i_j}$ is a Cartan-Serre generator corresponding to the simple root $\alpha_{i_j}$. Now consider the induced morphisms $\MM_{s_{i_{j}}\circ\lambda_{j-1},-h^\vee}\rightarrow \MM_{\lambda_{j-1},-h^\vee}$. If  the class $[f_{i_{1}}^{<\lambda,\alpha_{i_1}^\vee>+1}v_{\lambda}]\neq 0$, then thanks to $(iii)$,
 $[v_{\lambda_1}]$ is not annihilated by any non-zero element of $\cZ_-$, hence, as we have just argued,
$H^{\infty/2+0}_{DS}(L\fn_+,\MM_{\lambda_1, -h^\vee})$ is a free $\cZ_-$-module generated by
$[v_{\lambda_1}]$. Now an obvious inductive argument shows that it suffices to prove the
following:
\begin{equation}
\label{non-vanish-sing-vect}
f_i^n v_\mu\text{is not a coboundary.}
\end{equation}
{\em Proof.} Denote by $\fm_\chi$ the Lie subalgebra of $U(\fn_+[t^{-1}]t^{-1})$ obtained
by shifting $\fn_+[t^{-1}]t^{-1}$ by the Drinfeld-Sokolov character, see sect.\ref{The Drinfeld-Sokolov reduction Definition and main results},
Remark~\ref{rem-on-dr-sok-funct}. Namely, we define $\fm_\chi$ to be generated (as a Lie algebra) by $e^\alpha\otimes t^n$ if $n\leq -1$, $\alpha\in\Delta_+$ but not simple, and $e^\alpha\otimes t^{-1}+1$ if $\alpha\in\Delta_+$ is simple.
The definition of the Drinfeld-Sokolov differential implies $f_i^n v_\mu$ is a coboundary only
if $f_i^n v_\mu\in \fm_\chi \MM_{\mu,-h^\vee}$. We have an $\fm_\chi$-module decomposition $\MM_{\mu,-h^\vee}\stackrel{\sim}{\rightarrow} U(\fm_\chi)\otimes U(\fh[t^{-1}]t^{-1})\otimes U(\fn_-[t^{-1}])$. Upon this identification $f_i^n v_\mu\in 1\otimes1\otimes U(\fn_-[t^{-1}])$; therefore  it cannot belong in $\fm_\chi \MM_{\mu,-h^\vee}$. This concludes the proof of (\ref{non-vanish-sing-vect}) and hence of Lemma~\ref{reduction-one}. $\qed$

\begin{quotation}
{\em From now on $\lambda$ will be assumed to be regular, integral, dominant.}
\end{quotation}
\subsubsection{ }
\label{the-second-reduction}
Recall that $\MM_{\lambda/z}$ denotes the quotient $\MM_{\lambda,-h^\vee}/\cZ_-\MM_{\lambda,-h^\vee}$ and
let $N_\lambda$ be the (unique) maximal proper graded submodule of $\MM_{\lambda/z}$. We have a short exact sequence
\begin{equation}
\label{Nl-Mlz-Vlz}
0\rightarrow N_\lambda\rightarrow \MM_{\lambda/z}\rightarrow\VV_{\lambda/z}\rightarrow 0.
\end{equation}
Here is the 2nd reduction:
\begin{lem}
\label{2nd-reduction}
If
$H^{\infty/2+0}_{DS}(L\fn_+,N_{\lambda})=0$, then $H^{\infty/2+0}_{DS}(L\fn_+,\MM_{\lambda, -h^\vee})$ is a free $\cZ_-$-module generated by $[v_\lambda]$.
\end{lem}
{\em Proof.} Frenkel and Gaitsgory have computed \cite{FG3}
\[
H^{\infty/2+i}_{DS}(L\fn_+,\VV_{\lambda/z})=\left\{
\begin{array}{ll}
\CC[v_\lambda]&\text{ if }i=0\\
0&\text{ otherwise}.
\end{array}\right.
\]
Therefore, the long exact sequence of cohomology applied to (\ref{Nl-Mlz-Vlz}) gives
\[
0\rightarrow H^{\infty/2+0}_{DS}(L\fn_+,N_{\lambda})\rightarrow
H^{\infty/2+0}_{DS}(L\fn_+,\MM_{\lambda/z})\rightarrow\CC\rightarrow0.
\]
Hence the implication
\[
H^{\infty/2+0}_{DS}(L\fn_+,N_{\lambda})=0\Rightarrow H^{\infty/2+0}_{DS}(L\fn_+,\MM_{\lambda/z})=\CC=\CC[v_\lambda].
\]
It remains to prove the implication
\begin{equation}
\label{last-red-in-red-2}
H^{\infty/2+0}_{DS}(L\fn_+,\MM_{\lambda/z})=\CC[v_\lambda]\Rightarrow
H^{\infty/2+0}_{DS}(L\fn_+,\MM_{\lambda,-h^\vee})\stackrel{\sim}{\rightarrow}\cZ_-.
\end{equation}
We have already had a chance (Theorem~\ref{centre-on-verma-frenk}) to cite the Frenkel-Gaitsgory result that $\MM_{\lambda,-h^\vee}$ is a free $\cZ_-$-module.
The (increasing) degree filtration of $\cZ_-$ induces a filtration of $\MM_{\lambda,-h^\vee}$. Thus we obtain a spectral sequence
\[
\{E^{pq}_r\}\Rightarrow H^{\infty/2+p+q}_{DS}(L\fn_+,\MM_{\lambda,-h^\vee}).
\]
All its terms are $\cZ_-$-modules and, in particular,
\[
\oplus_p E^{p,-p}_1=\cZ_-\otimes H^{\infty/2+0}_{DS}(L\fn_+,\MM_{\lambda/z}).
\]
The L.H.S of (\ref{last-red-in-red-2}) implies
\[
\oplus_p E^{p,-p}_1\stackrel{\sim}{\rightarrow}\cZ_-.
\]
The higher terms of the spectral sequence are subquotients of the latter, but we already know, (\ref{sing-vect-cocycle}),  that $v_\lambda$ determines a cohomology class, which ``survives'' to the end, e.g. (\ref{non-vanish-sing-vect}); hence the higher terms of the spectral sequence are, in fact, quotients of $\cZ_-$.  The character formula (\ref{arakawa-coho-verma}) implies that all these quotients are isomorphic to $\cZ_-$. $\qed$

It remains to prove
\begin{prop}
\label{prop-to-prove-thm}
\[
H^{\infty/2+0}_{DS}(L\fn_+,N_{\lambda})=0.
\]
\end{prop}

\subsubsection{ }
\label{ds-hwv-nl-vanish}
The highest weight $\lambda$ being dominant, integral, regular, we obtain an exact sequence
\[
\bigoplus_{i=1}^{\text{rk}\fg}\MM_{s_i\circ\lambda,-h^\vee}\rightarrow N_\lambda.
\]
Let $N_i=\text{Im}\{\MM_{s_i\circ\lambda,-h^\vee}\rightarrow N_\lambda\}$. It is a highest weight module generated by $v_i$, the image of $v_{s_i\circ\lambda}$. Its twisted degree is $\langle s_i\circ\lambda,D-\rho^\vee\rangle$. To unburden the notation we will set $d_\mu=\langle\mu,D-\rho^\vee\rangle$.
Here is one approximation to Proposition~\ref{prop-to-prove-thm}.
\begin{prop}
\label{ds-of-hw-ni-vanish}
$H^{\infty/2+0}(L\fn_+,N_i)_{d_{s_i\circ\lambda}}=0$.
\end{prop}

The proof will given in sect. \ref{proof-proper-of-1st-approx} after some preparatory work in sect. \ref{weights-vs-ds}.

\subsubsection{ }
\label{weights-vs-ds}
Let $\LL_{\mu,-h^\vee}$ be the irreducible quotient of $\MM_{\mu,-h^\vee}$.
\begin{lem}
\label{twist-degr-comp-ser}
If $M\in \hat{\mathcal{O}}_k$, then
\[
H^{\infty/2+0}_{DS}(L\fn_+,M)=\bigoplus_{d\geq d_{min}}H^{\infty/2+0}_{DS}(L\fn_+,M)_d,
\]
where $d_{min}=\text{min}\{d_\mu:\;[M:\LL_{\mu,-h^\vee}]\neq0\}$.
\end{lem}
{\em Proof:} this assertion is an immediate consequence of the right exactness of the functor $H^{\infty/2+0}_{DS}(L\fn_+,?)$
(Theorem~\ref{dr-sok-as-derived}(3)) and the character formula (\ref{arakawa-coho-verma}). $\qed$

\begin{lem}
\label{lem-on-aff-weyl}
Let $\lambda$ be integral, regular, dominant.

(i) If $[\MM_{s_i\circ\lambda,-h^\vee}:\LL_{\mu,-h^\vee}]\neq0$ and $d_\mu\leq d_{s_i\circ\lambda}$, then $d_\mu=d_{s_i\circ\lambda}$ and
either $\mu=s_i\circ\lambda$ or $\mu=t_{\alpha_i^\vee\circ\lambda}$, where $t_{\alpha_i^\vee}=s_{-\alpha_i+\delta}s_{\alpha_i}$.

(ii) $t_{\alpha_i^\vee\circ\lambda}$ maximal among weights $\mu\neq s_i\circ\lambda$ such that $[\MM_{s_i\circ\lambda,-h^\vee}:\LL_{\mu,-h^\vee}]\neq0$.

(iii) $[\MM_{s_i\circ\lambda,-h^\vee}:\LL_{s_i\circ\lambda,-h^\vee}]=[\MM_{s_i\circ\lambda,-h^\vee}:\LL_{t_{\alpha_i^\vee}\circ\lambda-h^\vee}]=1$.
\end{lem}

{\em Proof of (i)} By virtue of the Kac-Kazhdan theorem \cite{KK}, $[\MM_{s_i\circ\lambda,-h^\vee}:\LL_{\mu,-h^\vee}]\neq0$ implies
$\mu\in W\circ\lambda -\ZZ_+\delta$. It is clear that $\langle w\circ\lambda,-\rho^\vee\rangle \geq \langle s_i\circ\lambda,-\rho^\vee\rangle$ unless $w=id$, and so
$d_\mu\geq d_{s_i\circ\lambda}$ unless $\mu=\lambda \text{ mod }\ZZ\delta$. Furthermore, $\langle w\circ\lambda,-\rho^\vee\rangle > \langle s_i\circ\lambda,-\rho^\vee\rangle$, hence $d_\mu>d_{s_i\circ\lambda}$, if $l(w)>1$.

If $\mu=s_j\circ\lambda-n\delta$, then $n\geq 1$, hence $d_\mu>d_{s_i\circ\lambda}$ unless $j=i$, $n=0$, in which case of course
$\mu=s_i\circ\lambda$.

Finally, if $\mu=\lambda$, then $\mu=s_{-\alpha_i+\delta}\circ s_i\circ\lambda$, which shows at once that  $d_\mu=d_{s_i\circ\lambda}$.

{\em Proof of (ii)}: this important for what follows  but technical assertion will be proved in sect.~\ref{proof-etchn-lemma(2)}.

{\em Proof of (iii)}: The equality $[\MM_{s_i\circ\lambda,-h^\vee}:\LL_{s_i\circ\lambda,-h^\vee}]=1$ is obvious. As to the equality
$[\MM_{s_i\circ\lambda,-h^\vee}:\LL_{t_{\alpha_i^\vee}\circ\lambda,-h^\vee}]=1$, item (ii) implies that the Shapovalov  form on $(\MM_{s_i\circ\lambda,-h^\vee})_{t_{\alpha_i^\vee}\circ\lambda}$ has a simple zero, and a standard application of the Jantzen filtration gives (iii).
$\qed$
\begin{lem}
\label{elimination-of-t-a-l}
If $[N_i:\LL_{\mu,-h^\vee}]\neq 0$ and $\mu\neq s_i\circ\lambda$, then $d_\mu> d_{s_i\circ\lambda}$.
\end{lem}
{\em Proof.} By  Lemma~\ref{lem-on-aff-weyl}(i), we have to eliminate the possibility of $\mu=t_{\alpha_i^\vee}\circ\lambda$.
Lemma~\ref{lem-on-aff-weyl}(iii) implies that the occurrence of $\LL_{\mu,-h^\vee}$ in the composition series of $N_i$ is due to the
composite morphism
\[
\MM_{t_{\alpha_i^\vee}\circ\lambda,-h^\vee}\hookrightarrow\MM_{s_i\circ\lambda,-h^\vee}\rightarrow N_i\subset\MM_{\lambda/z}
=\MM_{\lambda,-h^\vee}/\cZ_-\MM_{\lambda,-h^\vee}.
\]
But the latter lifts to a morphism $\MM_{t_{\alpha_i^\vee}\circ\lambda,-h^\vee}\rightarrow \MM_{\lambda,-h^\vee}$,  hence
defines an element of $\text{End}_{\widehat{\fg}}(\MM_{\lambda,-h^\vee})$. (This is because $\MM_{t_{\alpha_i^\vee}\circ\lambda,-h^\vee}$ and $\MM_{\lambda,-h^\vee}$ are isomorphic as $\widehat{\fg}$-modules.) By Theorem~\ref{centre-on-verma-frenk}, this element
equals an element of the center $\cZ_-$, hence the above composite morphism is 0. $\qed$
\subsubsection{ }
\label{proof-proper-of-1st-approx}
To the proof of Proposition~\ref{ds-of-hw-ni-vanish}.

Dualizing (\ref{Nl-Mlz-Vlz}), see sect.~\ref{duality-for-graded-modules} for the definition, we obtain the short exact sequence
\begin{equation}
\label{ex-seq-for-Nl-dual}
0\rightarrow \LL_{\lambda,-h^\vee}\rightarrow \MM_{\lambda/z}^c\rightarrow N_\lambda^c\rightarrow 0.
\end{equation}
We have $\MM_{\lambda/z}^c\stackrel{\sim}{\rightarrow}\WW_{\lambda/z}$ (see Lemma~\ref{what-dual-does-to-wak} below). Since
\[
H^{\infty/2+i}_{DS}(L\fn_+,\LL_{\lambda,-h^\vee})=
H^{\infty/2+i}_{DS}(L\fn_+,\WW_{\lambda/z})=\left\{\begin{array}{ll}
\CC&\text{ if }i=0\\
0&\text{ otherwise }
\end{array}
\right.,
\]
an application of the long exact cohomology sequence to (\ref{ex-seq-for-Nl-dual}) gives
\begin{equation}
\label{vanish-of-dual-Nl}
H^{\infty/2+\bullet}_{DS}(L\fn_+,N_{\lambda}^c)=0.
\end{equation}
$N_i$ being a submodule of $N_\lambda$, we get a surjection $N_\lambda^c\rightarrow N_i^c$,
hence, by Theorem~\ref{dr-sok-as-derived}(3) and (\ref{vanish-of-dual-Nl}),
\begin{equation}
\label{vanish-of-dual-Nl-i}
H^{\infty/2+\bullet}_{DS}(L\fn_+,N_i^c)=0.
\end{equation}
\begin{lem}
\label{surj-l-n}
(i) The surjection $N_i\rightarrow\LL_{s_i\circ\lambda,-h^\vee}$ induces an isomorphism\newline
$H^{\infty/2+0}_{DS}(L\fn_+,N_i)_{d_{s_i\circ\lambda}}\rightarrow H^{\infty/2+0}_{DS}(L\fn_+,\LL_{s_i\circ\lambda,-h^\vee})_{d_{s_i\circ\lambda}}$.

(ii) The embedding $\LL_{s_i\circ\lambda,-h^\vee}\rightarrow N_i^c$ induces an isomorphism\newline
$H^{\infty/2+0}_{DS}(L\fn_+,\LL_{s_i\circ\lambda,-h^\vee})_{d_{s_i\circ\lambda}}\rightarrow H^{\infty/2+0}_{DS}(L\fn_+,N_i^c)_{d_{s_i\circ\lambda}}$.
\end{lem}

{\em Proof.} Let $V$ be the kernel of $N_i\rightarrow\LL_{s_i\circ\lambda,-h^\vee}$. Obviously, $[V:\LL_{s_i\circ\lambda,-h^\vee}]=0$. If so, Lemma~\ref{elimination-of-t-a-l} says that $[V:\LL_{\mu,-h^\vee}]=[V^c:\LL_{\mu,-h^\vee}]\neq 0$ implies $d_\mu>d_{s_i\circ\lambda}$.
Now  Lemma~\ref{twist-degr-comp-ser} implies $H^{\infty/2+0}_{DS}(L\fn_+,V)_{d_{s_i\circ\lambda}}=0$. It remains to
use the long cohomology sequence (in case (ii) combined with the inequality
$\dim H^{\infty/2+0}_{DS}(L\fn_+,\LL_{s_i\circ\lambda,-h^\vee})_{d_{s_i\circ\lambda}}
\leq 1$, which is an obvious consequence of the right exactness and formula
(\ref{arakawa-coho-verma})). $\qed$

{\em Conclusion of the proof of Proposition~\ref{ds-of-hw-ni-vanish}}: Lemma~\ref{surj-l-n} implies the composition $N_i\rightarrow\LL_{s_i\circ\lambda,-h^\vee}\rightarrow N_i^c$
 induces an
isomorphism $H^{\infty/2+0}_{DS}(L\fn_+,N_i)_{d_{s_i\circ\lambda}}\rightarrow
H^{\infty/2+0}_{DS}(L\fn_+,N_i^c)_{d_{s_i\circ\lambda}}$. By (\ref{vanish-of-dual-Nl-i}),
the target of this map is 0, hence so is its source. $\qed$

\subsubsection{ } Recall, Remark~\ref{rem-on-dr-sok-funct}, that
$H^{\infty/2+\bullet}_{DS}(L\fn_+, ?)=H^{\infty/2+\bullet}(L\fn_+, ?\otimes\CC_\chi)$.
Here is another approximation to Proposition~\ref{prop-to-prove-thm}, one that replaces the semi-infinite
with ordinary Lie algebra homology.
\begin{prop}
\label{second--aprox-prop}
$H_0(\fn_+[t^{-1}]t^{-1},N_\lambda\otimes\CC_\chi)=0.$
\end{prop}
The proof will occupy the rest of this subsection with a proof of one auxiliary assertion spilling over to sect.~\ref{ Proof of whats-needed-spectr-conv-ha}.

Let $\{U_i(\widehat{\fg})\}$ be the standard filtration of the universal enveloping
algebra. The latter is acted upon by $\rho^\vee$, and we set
$U_i(\widehat{\fg})[j]=\{u\in U_i(\widehat{\fg}):\; [\rho^\vee,u]=-j\}$. The collection
of subspaces $K_n U_i(\widehat{\fg})=\sum_{i+j\leq n}U_i(\widehat{\fg})[j]$ defines a filtration of $U(\widehat{\fg})$ known as the {\em Kazhdan filtration}, cf. \cite{GG,A5}.
We have $\operatorname{gr}_K U(\widehat{\fg})=S(\widehat{\fg})$.

Any Lie subalgebra $\fa\subset \widehat{\fg}$ carries an induced filtration $\{K_n\fa\}$ so that $\operatorname{gr}_K\fa=S(\fa)$, and the embedding $\fa\hookrightarrow \widehat{\fg}$ induces
the embedding $S(\fa)\hookrightarrow S(\widehat{\fg})$.

Denote by $v_i$ a highest weight vector of $N_i$. Define
$K_n (N_\lambda\otimes\CC_\chi)=\sum_i(K_nU(\widehat{\fg}))(v_i\otimes 1)$; this is a Kazhdan filtration of $N_\lambda\otimes\CC_\chi$ compatible with that of $U(\fn_+[t^{-1}]t^{-1})$, and so
$\operatorname{gr}_K (N_\lambda\otimes\CC_\chi)$ is an  $S(\fn_+[t^{-1}]t^{-1})$-module generated by a collection $\{v_i\}$. Note that a similar definition applies to any finitely generated $\widehat{\fg}$-module.

This defines a filtration on the standard homology complex,
$\{K_nC_\bullet(\fn_+[t^{-1}]t^{-1},N_\lambda\otimes\CC_\chi)\}$. In the arising spectral sequence $\{E^r_{\bullet\bullet}\}$, the term $E^1_{\bullet\bullet}$ is the homology of the Koszul
complex $C_\bullet(\{e^\alpha\otimes t^n,\alpha\in\Delta_+,n<0\},\operatorname{gr}_{K}( N_\lambda\otimes\CC_\chi)$. Proposition~\ref{ds-of-hw-ni-vanish} implies
 (use the definition of the Drinfeld-Sokolov differential) that the symbol of each $v_i\otimes 1$ belongs to $ \fn_+[t^{-1}]t^{-1}\operatorname{gr}_{K} (N_\lambda\otimes\CC_\chi)$. Therefore
\[
\bigoplus_{p=-\infty}^{+\infty}E^1_{-p,p}=0.
\]
 Proposition~\ref{second--aprox-prop} will follow once   the convergence of the
 spectral sequence is proved. Since we are only interested in the vanishing of the zeroth homology group,  it suffices to prove that for each twisted degree $d$, sect.~\ref{Twisted grading and filtration dr-sok},
\begin{equation}
\label{whats-needed-spectr-conv}
H_0(K_pC_\bullet(\fn_+[t^{-1}]t^{-1},N_\lambda\otimes\CC_\chi)_d)=0\text { if }p\ll 0.
\end{equation}

\subsubsection{ Proof of (\ref{whats-needed-spectr-conv}):}
\label{ Proof of whats-needed-spectr-conv-ha}
 Set $M=\oplus_i\MM_{s_i\circ\lambda,-h^\vee}$. Consider the standard homology complex
$\{C_\bullet(\fn_+[t^{-1}]t^{-1},M\otimes\CC_\chi)\}$ and the Kazhdan filtration  $\{K_nC_\bullet(\fn_+[t^{-1}]t^{-1},N_\lambda\otimes\CC_\chi)\}$.

\begin{lem}
\label{vanish-homol-on-tw-verma-filtr}
For each twisted degree $d$
\[
H_\bullet(K_pC_\bullet(\fn_+[t^{-1}]t^{-1},M\otimes\CC_\chi)_d)=0\text { if }p\ll 0.
\]
\end{lem}
{\em Proof.} For the purpose of this proof let $C_\bullet(M)=C_\bullet(\fn_+[t^{-1}]t^{-1},M\otimes\CC_\chi)$. Consider a short exact sequence of complexes
\[
0\rightarrow K_pC_\bullet(M)_d\rightarrow C_\bullet(M)_d
\rightarrow C_\bullet(M)_d/K_pC_\bullet(M)_d\rightarrow 0.
\]
\begin{sloppypar}
Note that $M$ being free over $\fn_+[t^{-1}]t^{-1}$, $H_i(\operatorname{gr}_KC_\bullet(M))=0$ if $i\neq0$; furthermore, it is clear that $\dim H_0(\operatorname{gr}_KC_\bullet(M)_d)=\dim\operatorname{gr}_K H_0(C_\bullet(M)_d)<\infty$ for each $d$.
The Kazhdan filtration on $C_\bullet(M)_d/K_pC_\bullet(M)_d$ being regular,
$H_0(\operatorname{gr}_K C_\bullet(M)_d/K_pC_\bullet(M)_d)=\operatorname{gr}_K H_0( C_\bullet(M)_d/K_pC_\bullet(M)_d)$. But $\operatorname{gr}_K C_\bullet(M)_d/K_pC_\bullet(M)_d$ id a direct summand of the complex $C_\bullet(M)_d$. The finite dimensionality, $\dim H_0(\operatorname{gr}_KC_\bullet(M)_d)_d<\infty$, implies that the rightmost map in the exact sequence above, $C_\bullet(M)_d
\rightarrow C_\bullet(M)_d/K_pC_\bullet(M)_d$ induces an isomorphism of homology  $H_\bullet(C_\bullet(M)_d)
\rightarrow H_\bullet(C_\bullet(M)_d/K_pC_\bullet(M)_d)$ for all $p\ll 0$. Lemma~\ref{vanish-homol-on-tw-verma-filtr} follows. $\qed$
\end{sloppypar}

\begin{sloppypar}
Now an obvious surjection $M\rightarrow N_\lambda$ induces a surjection of complexes $K_pC_\bullet(\fn_+[t^{-1}]t^{-1},M\otimes\CC_\chi)_d\rightarrow K_pC_\bullet(\fn_+[t^{-1}]t^{-1},N_\lambda\otimes\CC_\chi)_d)$. Since the homology functor is exact, we obtain a surjection
$H_0(K_pC_\bullet(\fn_+[t^{-1}]t^{-1},M\otimes\CC_\chi)_d) \rightarrow H_0(K_pC_\bullet(\fn_+[t^{-1}]t^{-1},N_\lambda\otimes\CC_\chi)_d)$.
Thus Lemma~\ref{vanish-homol-on-tw-verma-filtr} implies (\ref{whats-needed-spectr-conv}).       The proof of Proposition~\ref{second--aprox-prop}
is now complete.
\end{sloppypar}
\subsubsection{ }
\label{relative-coho}
We will now undertake a short digression on relative cohomology. For any $M\in\cO_k$, consider
\begin{multline}
\nonumber
C^{\infty/2+\bullet}(L\fn_+,\fn_+[t^{-1}]t^{-1}, M\otimes\CC_\chi)\\
=\frac{C^\bullet(L\fn_+, M\otimes\CC_\chi)}{\text{span}\{(\phi_{\alpha})_{(-n)}c,[d_{st},(\phi_{\alpha})_{(-n)}]c,\;
\alpha\in\Delta_+,n>0,c\in C^\bullet(L\fn_+, M\otimes\CC_\chi)\}}.
\end{multline}
It is a quotient complex of the Drinfeld-Sokolov complex $(C^{\infty/2+\bullet}(L\fn_+, M\otimes\CC_\chi),d_{st})$, cf. Remark~\ref{rem-on-dr-sok-funct},
known as the {\em relative semi-infinite cohomology complex}
of $L\fn_+$ with coefficients in $M\otimes\CC_\chi)$. The corresponding relative cohomology groups will be denoted by $H^{\infty/2+i}(L\fn_+,\fn_+[t^{-1}]t^{-1}, M\otimes\CC_\chi)$; note that $i\geq 0$.
\begin{thm}
\label{iso-ds-funct-relat-funct}
The functors $H^{\infty/2+0}(L\fn_+,\fn_+[t^{-1}]t^{-1}, ?\otimes\CC_\chi)$ and
$H^{\infty/2+0}_{DS}(L\fn_+, ?)$ are naturally isomorphic.
\end{thm}
The rest of this subsection will be devoted to the proof of this theorem.

\begin{lem}
\label{coinc-on-vermas}
If $P$ is a module with a Verma filtration, then
\[
H^{\infty/2+\bullet}(L\fn_+,\fn_+[t^{-1}]t^{-1}, P\otimes\CC_\chi)\stackrel{\sim}{\rightarrow}H^{\infty/2+\bullet}_{DS}(L\fn_+, P).
\]
\end{lem}
\begin{sloppypar}
{\em Proof.} Consider the Serre-Hochschild spectral sequence \cite{V} $\{E_r^{pq}\}$ with
$E_1^{pq}=H_{-p}(\fn_+[t^{-1}]t^{-1},P\otimes\CC_\chi\otimes\Lambda^q\fn_+[t]^*)$.
Since $P$ is filtered by Verma modules, it is free over $\fn_+[t^{-1}]t^{-1}$, hence
$E_1^{pq}=0$ unless $p=0$, in which case it equals $H_{0}(\fn_+[t^{-1}]t^{-1},P\otimes\CC_\chi\otimes\Lambda^\bullet\fn_+[t]^*)$, and
 the spectral sequence collapses. In fact, a little more is valid.

Note that although the differential does not
preserve conformal grading (this is $\CC_\chi$'s doing), it preserves the increasing conformal filtration
$\{F_nC^{\infty/2+p+q}_{DS}(L\fn_+, P)\}$, cf. sect.~\ref{Gradings and character formulas}.
It is easy to see that the aforementioned vanishing takes place upon restriction to each term
of this filtration: $(F_nE)^{pq}_1=0$ unless $p=0$. Now notice that  the filtration that leads to
$\{E_r^{pq}\}$ induces a finite filtration of  $F_nC^{\infty/2+p+q}_{DS}(L\fn_+, P)$ for each $n$. A quick diagram chase then shows the convergence of the spectral sequence: $\{E_r^{pq}\}\Rightarrow H^{\infty/2+p+q}_{DS}(L\fn_+, P)$.

Since the spectral sequence collapses, we obtain that the natural map
$H^{\infty/2+\bullet}_{DS}(L\fn_+, P)\rightarrow H^\bullet(E^1_{0\bullet})$ defines an isomorphism
$H^{\infty/2+\bullet}_{DS}(L\fn_+, P)\stackrel{\sim}{\rightarrow}H^\bullet(E^1_{0\bullet})$
\end{sloppypar}
On the other hand, by definition, the complexes $E^1_{0\bullet}$ and $(C^\bullet(L\fn_+, \fn_+[t^{-1}]t^{-1},M\otimes\CC_\chi),d_{st})$ are identical. $\qed$

\bigskip
Since the Drinfeld-Sokolov functor $H^{\infty/2+0}_{DS}(L\fn_+, ?)$ is right exact and the class of modules with a Verma filtration is adapted, it remains to prove that
the functor $H^{\infty/2+0}(L\fn_+,\fn_+[t^{-1}]t^{-1}, ?\otimes\CC_\chi)$ is right exact.

Instrumental for proving the right exactness is Proposition 3.7.1 of  \cite{A5}. This result is but a jet scheme version of [GG], sect.5 and 6, and it deals with the following
situation. Let $\fb\subset\fa$ is a pair ``Lie algebra, Lie subalgebra'', and let a vector
space $M$ be acted upon by $S(\fa)$ as a commutative algebra, and by $\fb$ as a Lie algebra. Denote by $a\bullet ?$ and $\{b,?\}$ the respective actions by $a\in\fa$ and $b\in\fb$. We shall say that the 2 structures are compatible if $\{b,a\bullet m\}=
[b,a]\bullet m+a\bullet\{b,m\}$.
\begin{lem} (cf. Proposition 3.7.1 of  \cite{A5}.)
\label{compil-arak-result}
 Let $M$ be acted upon by $S(\widehat{\fg}/\fg[t])$ as a commutative algebra and by $\fg[t]$ as a Lie algebra, and let the 2 structures  be compatible. If the shifted subalgebra
 $S(\fm_\chi)=S(\text{span}\{e_\alpha\otimes t^n+\chi(e_\alpha\otimes t^n),\;n<0,\alpha\in\Delta_+\})$ (cf. Remark~\ref{rem-on-dr-sok-funct}) satisfies $S(\fm_\chi)^N M=0$ for some $N>0$, then
 \[
 H^{i}(\fn_+[t],M)=\left\{\begin{array}{ll}M^{\fn_+[t]}&\text{ if }i=0\\
 0&\text{ otherwise.}
 \end{array}
  \right.
 \]
 \end{lem}

This lemma will be applied via the decreasing Li filtration \cite{Li2}. Given  a vertex algebra $V$ and a $V$-module $M$, Li defines $\operatorname{L}_nM$ to be the linear span
$x^1_{(-i_1-1)}\cdots x^l_{(-i_l-1)}m$ with  $x^j\in V$, $1\leq j\leq l$, $m\in M$, and
$\sum_ji_j\geq n$. Li checks that the graded objects satisfy: $\operatorname{gr}_LV$
is a vertex Poisson algebra, and
$\operatorname{gr}_LM$ is a $\operatorname{gr}_LV$-module.
Whatever else
this means (see also \cite{FBZ} for details), this implies that
if $M\in\hat{\cO}_k$, then $\operatorname{gr}_LM$ is an $S(\widehat{\fg}/\fg[t])$-module
and a $\fg[t]$-module, and the 2 structures are compatible.

\begin{cor}
\label{vanish-rel-coho} For each $M\in\hat{\cO}_k$,
$H^{\infty/2+i}(L\fn_+,\fn_+[t^{-1}]t^{-1}, M\otimes\CC_\chi)=0$ if $i\neq 0$. Furthermore, $H^{\infty/2+0}(L\fn_+,\fn_+[t^{-1}]t^{-1}, M\otimes\CC_\chi)=0$ is filtered, and $\operatorname{gr}H^{\infty/2+0}(L\fn_+,\fn_+[t^{-1}]t^{-1}, M\otimes\CC_\chi)=(\operatorname{gr}M/\fm_\chi\operatorname{gr}M)^{\fn_+[t]}$
\end{cor}

{\em Proof of Corollary~\ref{vanish-rel-coho}.} Consider the decreasing Li \cite{Li2} filtration of the relative semi-infinite cohomology complex. A spectral sequence arises, and its second term
is $H^{\bullet}(\fn_+[t],\operatorname{gr}_LM/\fm_\chi\operatorname{gr}_LM)$. As explained above, this brings us into the situation of Lemma~\ref{compil-arak-result}. An application of Lemma~\ref{compil-arak-result} concludes the proof.$\qed$

\bigskip
The right exactness of the functor $H^{\infty/2+0}(L\fn_+,\fn_+[t^{-1}]t^{-1}, ?\otimes\CC_\chi)$ is now easy. Given an exact sequence
\[
0\rightarrow K\rightarrow M\rightarrow N\rightarrow 0,
\]
we consider its Li  graded version (as in the proof above)
\[
0\rightarrow \operatorname{gr}_LK\rightarrow \operatorname{gr}_LM\rightarrow \operatorname{gr}_LN\rightarrow 0,
\]
 obtain an exact sequence of coinvariants (0th Lie algebra homology)
\[
 \operatorname{gr}_LK/\fm_\chi \operatorname{gr}_LK\rightarrow
 \operatorname{gr}_LM/\fm_\chi \operatorname{gr}_LM\rightarrow
 \operatorname{gr}_LN/\fm_\chi \operatorname{gr}N_L\rightarrow 0,
\]
and letting $U=\text{Ker}\{\operatorname{gr}_LK/\fm_\chi \operatorname{gr}_LK\rightarrow
 \operatorname{gr}_LM/\fm_\chi \operatorname{gr}_LM\}$, an exact cohomology sequence
\[
 (\operatorname{gr}_LM/\fm_\chi \operatorname{gr}_LM)^{\fn_+[t]}\rightarrow
 (\operatorname{gr}_LN/\fm_\chi \operatorname{gr}_LN)^{\fn_+[t]}\rightarrow
 H^{1}(\fn_+[t],U)\stackrel{\text{Lemma}\ref{compil-arak-result}}{=}0.
\]
\begin{sloppypar}
By virtue of Corollary~\ref{vanish-rel-coho}, this shows that the graded functor, $\operatorname{gr}H^{\infty/2+0}(L\fn_+,\fn_+[t^{-1}]t^{-1}, ?\otimes\CC_\chi)$, is right exact. Hence so is the functor itself. This concludes the proof of Theorem~\ref{iso-ds-funct-relat-funct}. $\qed$
\end{sloppypar}
\subsubsection{ }
\label{thm-from-2nd-approx}
We can now derive Proposition~\ref{prop-to-prove-thm} from its ``2nd approximation'', Proposition~\ref{second--aprox-prop}. Recall that we want
to show $H^{\infty/2+0}_{DS}(L\fn_+, N_\lambda)=0$. According to Theorem~\ref{iso-ds-funct-relat-funct}, this is equivalent to
$H^{\infty/2+0}(L\fn_+,\fn_+[t^{-1}]t^{-1}, N_\lambda\otimes\CC_\chi)=0$. But we have the vanishing already at the level of the complex: indeed, by the definition
of the relative semi-infinite cohomology, sect.~\ref{relative-coho},
$C^{\infty/2+0}(L\fn_+,\fn_+[t^{-1}]t^{-1}, N_\lambda\otimes \CC_\chi)=N_\lambda/\fm_\chi N_\lambda$, which is 0 thanks to Proposition~\ref{second--aprox-prop}. This concludes the proof of Proposition~\ref{prop-to-prove-thm}, hence of Theorem~\ref{freeness-ds-for-verma} -- modulo, that is, Lemma~\ref{lem-on-aff-weyl}(ii). Its proof is as follows:
\subsubsection{ }
\label{proof-etchn-lemma(2)}
We will be freely using the results and notation from \cite{Lus}. Let $A_w$ be the alcove attached to an affine Weyl group element $w\in\widehat{W}$;
for example $A_w=w^{-1}A_1$. Choose quarters $C^+_v$ (cf. \cite{Lus}, sect. 1.1) as the translates of the fundamental chamber of $\fg$.

Defined in \cite{Lus}, sect. 1.4, is a number, $d(A,B)$, for any two alcoves $A$, $B$.
\begin{lem}
\label{lem-lust-1}
$d(A_w,A_1)=l^{\frac{\infty}{2}}(w)$, where
\[
l^{\frac{\infty}{2}}(w)=\#\{\alpha\in\hat{\Delta}^{re}_+\cap w^{-1}(\hat{\Delta}^{re}_-),\;\bar{\alpha}\in\Delta_+\}-
\#\{\alpha\in\hat{\Delta}^{re}_+\cap w^{-1}(\hat{\Delta}^{re}_-),\;\bar{\alpha}\in\Delta_-\},
\]
and $\hat{\fh}^*\ni\lambda\mapsto\bar{\lambda}\in\fh$ is the restriction map, and $\hat{\Delta}_{\pm}^{re}$ is the set of positive (negative) real roots of $\widehat{\fg}$.
\end{lem}
{\em Proof.} Use \cite{Lus}, sect.~1.4 to proceed by induction on the length $l(w)$. If $l(w)=1$, the assertion follows from the definition.
Let $w=ys_i$ with $l(w)=l(y)+1$. Then $A_y$ and $A_w$ are next to each other and are separated by the hyperplane $y^{-1}(H_i)$. Hence
$d(A_w, A_1)=d(A_y, A_1)\pm1$. By definition, $d(A_w, A_1)=d(A_y, A_1)+1$ iff $\overline{y^{-1}(\alpha_i)}\in\Delta_+$. Since
$\hat{\Delta}^{re}_+\cap w^{-1}(\hat{\Delta}^{re}_-)=\hat{\Delta}^{re}_+\cap y^{-1}(\hat{\Delta}^{re}_-)\sqcup\{y^{-1}(\alpha_i)\}$,
the assertion follows. $\qed$
\begin{lem}
\label{lem-lust-2}
Let $\lambda\in\hat{\fh}^*$ be a weight at the critical level such that $\bar{\lambda}$ is dominant regular, $w\in\widehat{W}$, $\alpha\in\hat{\Delta}^{re}_+$. If there is an embedding $\MM_{s_\alpha w\circ\lambda}\rightarrow\MM_{w\circ\lambda}$, then
$l^{\frac{\infty}{2}}(s_\alpha w)>l^{\frac{\infty}{2}}(w)$.
\end{lem}
{\em Proof.} By \cite{Lus}, (1.4.1), $l^{\frac{\infty}{2}}(s_\alpha w)=l^{\frac{\infty}{2}}(w)+d(A_{s_\alpha w}, A_w)$. Note that
$A_{s_\alpha w}=A_w\sigma_{H_{w^{-1}(\alpha)}}$. By assumption, $\overline{w^{-1}(\alpha)}\in\Delta_+$. Therefore, by \cite{Lus}, Lemma~3.1,
$d(A_{s_\alpha w}, A_w)>0$. $\qed$

\bigskip
Now Lemma~\ref{lem-on-aff-weyl}(ii) follows from
\begin{lem}
\label{lem-lust-3}
For $i\neq0$ (equivalently, $\alpha_i$ is a simple root of the finite root system), $l^{\frac{\infty}{2}}(s_i)=1$ and
$l^{\frac{\infty}{2}}(t_{\alpha_i^\vee})=2$.
\end{lem}
{\em Proof.} The first assertion is obvious. The second is also elementary: each element of $\beta\in\Delta_+$ will contribute
$\langle\beta,\alpha_i^\vee\rangle$ toward $l^{\frac{\infty}{2}}(t_{\alpha_i^\vee})$; overall we obtain $2\langle\rho,\alpha_i^\vee\rangle$,
which is 2. $\qed$

The proof of Theorem~\ref{freeness-ds-for-verma} is now complete.

\section{Proof of Theorem~\ref{main-theorem-intro}}
\label{Proof of Theorem-main-theorem-intro-homo}
\subsection{Resolution}
\label{resolution} We shall work in the setting of sect.~\ref{examples-of-d-ch-tw-mod}

The Cousin resolution of $\cL_{\nu_0}$ w.r.t. the filtration of $X$ by $\{X_w,\,w\in W\}$ reads
\begin{equation}
\label{cousin} 0\rightarrow \cL_{\nu_0}\rightarrow i_{{\rm id},+}i_{\rm id}^*\cL_{\nu_0}
\rightarrow \oplus_{w\in W^{(1)}}i_{w,+}i_w^*\cL_{\nu_0} \rightarrow \oplus_{w\in
W^{(2)}}i_{w,+}i_w^*\cL_{\nu_0} \rightarrow\cdots.
\end{equation}

Note that
\begin{equation}
\label{from-local-to-contr-verma}
 \Gamma(X, i_{w,+}i_w^*\cL_{\nu_0})=M^{c}_{w\circ\nu_0},
 \end{equation}
 the contragredient Verma module.

 Let $\cM^w_{\nu(z)}$ stand for $\cZ hu_{\nu(z)}(i_{w,+}i_w^*\cL_{\nu_0})$, see (\ref{left-adj-to-forg-functor}).

The functor $\cZ hu_{\nu(z)}$ is exact \cite{AChM}, and its application to  (\ref{cousin})
gives a resolution in the category of $\cD^{ch,tw}_X$-modules
\begin{equation}
\label{cousin-groth-chiral} 0\rightarrow
\cL_{\nu(z)}^{ch}\rightarrow \cM^{{\rm id}}_{\nu(z)} \rightarrow
\oplus_{w\in W^{(1)}}\cM^{w}_{\nu(z)} \rightarrow \oplus_{w\in
W^{(2)}}\cM^w_{\nu(z)} \rightarrow\cdots.
\end{equation}
Lemma~\ref{vanish-higher-delta-fnct} implies
\[
H^i(X,\cM^w_{\nu(z)})=0\text{ if }i>0.
\]
Since the class of sheaves with vanishing higher cohomology is
adapted to the functor of global sections,  the complex
\begin{equation}
\label{cousin-groth-chiral-glob-sect} 0\rightarrow \Gamma(X,
\cM^{{\rm id}}_{\nu(z)}) \rightarrow \oplus_{w\in
W^{(1)}}\Gamma(X,\cM^w_{\nu(z)}) \rightarrow \oplus_{w\in
W^{(2)}}\Gamma(X,\cM^w_{\nu(z)}) \rightarrow\cdots.
\end{equation}
computes the cohomology $H^i(X,\cL_{\nu(z)}^{ch})$.

Now recall (\ref{def-wak-crit-lev-arb-centr-charge-tw-w}) that $\Gamma(X,\cM^w_{\nu(z)})$ is
precisely the $w$-twisted Wakimoto module $\WW^{w}_{\nu(z)}$. We can summarize our discussion
as follows.
\begin{lem}
\label{reform-wak-modul-proof} The cohomology of the complex
\begin{equation}
\label{reform-wak-modul-proof-compl} 0\rightarrow \WW^{{\rm
id}}_{\nu(z)} \rightarrow \oplus_{w\in W^{(1)}}\WW^w_{\nu(z)}
\rightarrow \oplus_{w\in W^{(2)}}\WW^w_{\nu(z)} \rightarrow\cdots.
\end{equation}
is isomorphic to $H^\bullet(X,\cL_{\nu(z)}^{ch})$.
\end{lem}
Both  items of Theorem~\ref{main-theorem-intro} follow from this
lemma, the first easily, the second after some work is done.

\subsection{ Proof of Theorem~\ref{main-theorem-intro}(1) and formula (\ref{homo=vers-of=thm-1.1.2}).}
\label{Proof-of-Theoremrefmain-theorem-intro(1)}  First of all, the
definition used in Theorem~\ref{main-theorem-intro}(1),
\[
\chi(\cL_{\nu(z)}^{ch})=\sum_{i=0}^{{\rm dim}X}(-1)^i{\rm ch}
H^i(X,\cL_{\nu(z)}^{ch}),
\]
has been made sense of:  $\chi(\cL_{\nu(z)}^{ch})$ being a sheaf of
 filtered, see (\ref{constr-of-zhu-filtr}), $\ghat$-modules, see Theorem~\ref{exist-tw-ch-do-flag} (4),
  the cohomology groups
  $H^i(X,\cL_{\nu(z)}^{ch})$ are objects of $\hat{\cO}_{\nu(z)\circ\pi}$, the full subcategory
  of $\hat{\cO}_{crit}$, cf. sect.~\ref{def-of-o-cat}; their
  formal characters, ${\rm ch}
H^i(X,\cL_{\nu(z)}^{ch})$, are discussed in sect.~\ref{Gradings and
character formulas}

 Lemma~\ref{reform-wak-modul-proof} implies
\begin{equation}
\label{char-weyl-1} \chi(\cL_{\nu(z)}^{ch})=\sum_{w\in
W}(-1)^{l(w)}{\rm ch}\WW^w_{\nu(z)}.
\end{equation}
Since, see (\ref{char-of-wak-mod-restr}),
\[
{\rm
ch}\WW^w_{\nu(z)}=e^{w\circ\nu_0}\prod_{\alpha\in\Delta_+}(1-e^{-\alpha})^{-1}\times
\prod_{n=1}^{+\infty}\prod_{\alpha\in\Delta_+}(1-e^{-\alpha-n\delta})^{-1}\times
\prod_{\alpha\in\Delta_-}(1-e^{-\alpha-n\delta})^{-1}
\]
(\ref{char-weyl-1}) can be rewritten as follows
\begin{align}
\label{char-weyl-2} \chi(\cL_{\nu(z)}^{ch})=&\left(\sum_{w\in
W}(-1)^{l(w)}e^{w\circ\nu_0}\prod_{\alpha\in\Delta_+}(1-e^{-\alpha})^{-1}\right)\\
\times&
\prod_{n=1}^{+\infty}\prod_{\alpha\in\Delta_+}(1-e^{-\alpha-n\delta})^{-1}
\prod_{\alpha\in\Delta_-}(1-e^{-\alpha-n\delta})^{-1},\nonumber
\end{align}
which is  Theorem~\ref{main-theorem-intro}(1) in a slightly expanded
form. In order to obtain (\ref{homo=vers-of=thm-1.1.2}), we have to
let $e^\alpha\rightarrow 1$, $\alpha\in\Delta$, and set
$e^{-\delta}=q$. In the limit, the first factor in
(\ref{char-weyl-2}) equals ${\rm dim}V_{\nu_0}$ (the Weyl character
formula !), and (\ref{char-weyl-2}) becomes
\[
\chi_P(\cL_{\nu(z)}^{ch})={\rm
dim}V_{\nu_0}\prod_{n=1}^{+\infty}(1-q^n)^{-2{\rm dim}X},
\]
as desired.

\subsection{ Proof of Theorem~\ref{main-theorem-intro}(2): the case of a homogeneous character.}
\label{ Proof-of-Theorem-main-theorem-intro(2)} Since each
$\cL^{ch}_{\nu(z)}$ is $G$-equivariant,
$H^i(X,\cL^{ch}_{\nu(z)})\in\hat{\cO}_{{\rm
crit},\nu(z)\circ\pi}^{G}$, the full subcategory of
$\hat{\cO}_{crit}$ that consists of those $\ghat$-modules where the
action of $\fg\subset\ghat$ integrates to an action of $G$. This
category being semi-simple with a unique simple object
$\VV_{\nu(z)}$ \cite{FG3}, we obtain
\begin{equation}
\label{what-coho-looks-apriori}
H^i(X,\cL^{ch}_{\nu(z)})=\oplus_{j=1}^{m_i}\VV_{\nu(z)}[n_{ij}],
\end{equation}
for some nonnegative integers $m_i$ and $\{n_{ij}\}$, which we have
to determine. The meaning of $\VV_{\nu(z)}[n_{ij}]$ is as follows:
the LHS of (\ref{what-coho-looks-apriori}) is filtered via
(\ref{constr-of-zhu-filtr}); each $\VV_{\nu(z)}$ appearing on the
right inherits a filtration, and the inherited filtration can be
different from the {\it natural } one, see sect.~\ref{Gradings and
character formulas} by a shift; this shift is denoted by $[n_{ij}]$.

It is proved in \cite{FG3} that  the Drinfeld-Sokolov reduction
functor extracts these numbers:
\begin{equation}
\label{what-coho-looks-apriori-via-dr-sok}
H^{\infty/2+k}_{DS}\left(L\fn_+,H^i(X,\cL^{ch}_{\nu(z)})\right)=\left\{
\begin{array}{ll}
\oplus_{j=1}^{m_i}\CC[n_{ij}]&\text{ if }k=0\\
0&\text{ otherwise}.
\end{array}\right.
\end{equation}
It is easy to interpret (\ref{what-coho-looks-apriori-via-dr-sok})
is a computation of the cohomology of a certain double complex. Let
\[
\cC^{j}=\oplus_{w\in W^{(j)}}\WW^w_{\nu(z)},
\]
and extend (\ref{reform-wak-modul-proof-compl}) to a double complex
\begin{equation}
\label{double-complex}
K^{\bullet\bullet}=\oplus_{pq}K^{pq},\;K^{pq}=(C^{\infty/2+p}(L\fn_+,\cC^q),
d_{DS}+d),
\end{equation}
where by $d$ we have denoted the differential of complex
(\ref{what-coho-looks-apriori-via-dr-sok}).

Either of the two spectral sequences associated with $K^{pq}$
converges to  $H^{p+q}_{d_{DS}+d}(K^{\bullet\bullet})$ because
complex (\ref{what-coho-looks-apriori-via-dr-sok}) is of finite
length. What (\ref{what-coho-looks-apriori-via-dr-sok}) says is that
one of these spectral sequences collapses in the second term and
\begin{equation}
\label{coho-via-first-proof-of-2}
H^{i}_{d_{DS}+d}(K^{\bullet\bullet})=\oplus_{j=1}^{m_i}\CC[n_{ij}].
\end{equation}
We will compute another spectral sequence (not without a twist)
thereby proving the following
\begin{lem}
\label{another-spectr-seq-Kpq} If $\nu(z)=\nu_0/z$,
\[
H^{i}_{d_{DS}+d}(K^{\bullet\bullet})=\oplus_{w\in
W^{(i)}}\CC[\langle\nu_0-w\circ\nu_0,\rho^{\vee}\rangle].
\]
\end{lem}
This lemma along with
(\ref{what-coho-looks-apriori},\ref{what-coho-looks-apriori-via-dr-sok},\ref{coho-via-first-proof-of-2})
implies Theorem~\ref{main-theorem-intro}(2).

\subsubsection{Proof of Lemma~\ref{another-spectr-seq-Kpq}.}
\label{proof-of-lemma-anoth-spectr-seq} Since $\nu(z)=\nu_0/z$, each
$\cC^j$ is a graded $\ghat$-module, and we can dualize \footnote{this
the only place, where the homogeneity assumption $\nu(z)=\nu_0/z$ is
used}, see sect.~\ref{duality-for-graded-modules}, to obtain $(\cC^q)^c$ and
\[
\tilde{K}^{\bullet\bullet}=\oplus_{pq}\tilde{K}^{pq},\;K^{pq}=(C^{\infty/2+p}(L\fn_+,(\cC^q)^c),
d_{DS}+d).
\]
Since $\VV_{\nu(z)}$ is irreducible,
$(\VV_{\nu_0/z})^c=\VV_{\nu_0/z}$ and the argument that led above to
(\ref{coho-via-first-proof-of-2}) can be repeated to give us
\begin{equation}
\label{coho-via-first-proof-of-2-dual}
H^{i}_{d_{DS}+d}(\tilde{K}^{\bullet\bullet})=\oplus_{j=1}^{m_i}\CC[n_{ij}]=
H^{i}_{d_{DS}+d}(K^{\bullet\bullet}).
\end{equation}
Therefore, it suffices to prove Lemma~\ref{another-spectr-seq-Kpq}
with $K^{\bullet\bullet}$ replaced with
$\tilde{K}^{\bullet\bullet}$. In order to do this, consider that
spectral sequence $\{E^{pq}_r\}$ where
\[
E^{pq}_1=H^{\infty/2+q}_{DS}((\cC^p)^c).
\]
Now we wish to compute $\{E^{pq}_1\}$.
\begin{lem}
\label{what-dual-does-to-wak}
\[\WW^w_{\nu_0/z}\stackrel{\sim}{=}\MM_{w\circ\nu_0/z}^c.
\]
\end{lem}
Lemma~\ref{what-dual-does-to-wak} implies
\[
(\WW^w_{\nu_0/z})^c\stackrel{\sim}{=}\MM_{w\circ\nu_0/z}.
\]
Since $(\cC^p)^c=\oplus_{w\in W^{(p)}}(\WW^w_{\nu_0/z})^c$, now
known to be $\oplus_{w\in W^{(p)}}\MM_{w\circ\nu_0/z}$,
Corollary~\ref{dr-sok-of-restr-verma} gives
\begin{equation}
\label{1st-term-via-another-spectr-seq}
 E^{pq}_1=\left\{
\begin{array}{ll}\oplus_{w\in
W^{(i)}}\CC[\langle\nu_0-w\circ\nu_0,\rho^{\vee}\rangle]&\text{ if
}q=0\\
0&\text{ otherwise,} \end{array}\right.
\end{equation}
where $\CC[\langle\nu_0-w\circ\nu_0,\rho^{\vee}\rangle]$ is spanned
by the class of the highest weight vector of $\MM_{w\circ\nu_0/z}$;
note that the grading shift is precisely the one obtained by placing
a highest weight vector of $\MM_{\nu_0/z}$ in degree 0 component and
then using the twisted grading as defined in sect.~\ref{Twisted
grading and filtration dr-sok}. The following dimensional argument
shows that all higher differentials vanish. The cohomology classes
recorded in (\ref{1st-term-via-another-spectr-seq}) are represented
by the highest weight vectors ${\bf 1}_{w\circ\nu_0}\in
\WW^w_{\nu_0/z}$. The twisted degree, the one that is preserved by
the differential of the double complex (\ref{double-complex}), of
${\bf 1}_{w\circ\nu_0}$ is $\langle w\circ\nu_0,\rho^{\vee}\rangle$.
The differential $d_1$ that operates on $E^{\bullet\bullet}_1$ is
induced by the differential of  complex
(\ref{reform-wak-modul-proof-compl}). The latter (by the
construction of the Cousin resolution (\ref{cousin})) is a direct
sum of the maps
\[
\WW^w_{\nu_0/z}\rightarrow \WW^v_{\nu_0/z}\text{ with }v>w,
l(v)=l(w)+1.
\]
 Since $\nu_0$ is regular dominant,
$\langle w\circ\nu_0,\rho^{\vee}\rangle<\langle
v\circ\nu_0,\rho^{\vee}\rangle$ provided $w<v$. Hence $E^{p0}_1$ and
$E^{p-1 0}_1$ have different twisted degrees, which makes $d_1$ be
equal to zero. The vanishing of higher differentials is obvious.
This concludes the proof of Lemma~\ref{another-spectr-seq-Kpq}.

\subsubsection{Proof of Lemma~\ref{what-dual-does-to-wak}.}
\label{Proof of Lemma what-dual-does-to-wak } This proof is but a
version of the argument in \cite{F3}, sect. 9.5.2, just as
Lemma~\ref{what-dual-does-to-wak} is a  generalization of
Proposition~9.5.1 proved there.

 The universal property of the Verma
module implies the existence of a non-zero morphism
\[\MM_{w\circ\nu_0/z}\rightarrow (\WW^w_{\nu_0/z})^c.\]
Dualizing we obtain
\[\WW^w_{\nu_0/z}\rightarrow \MM_{w\circ\nu_0/z}^c.\]
In order to prove that this map is an isomorphism, it suffices to
prove that $\WW^w_{\nu_0/z}$ has a unique up to proportionality
singular vector; in which case it may only belong to $\CC{\bf
1}_{w\circ\nu_0}$. Formally, we want to prove that
\begin{equation}
\label{what we want-lemma-dual-does-to-wak}
(\WW^w_{\nu_0/z})^{\hat{\fn}_+}=\CC{\bf 1}_{w\circ\nu_0}.
\end{equation}
Consider the Lie subalgebra $Lw(\fn_+)\cap\hat{\fn}_+$. By
construction, $\WW^w_{\nu_0/z}$ is co-free as an
$Lw(\fn_+)\cap\hat{\fn}_+$-module; this sort of assertion has been
the cornerstone of the Wakimoto module theory since its inception,
cf. \cite{F3}, sect. 9.5.2. Informally speaking,
if $U^w$ is the maximal unipotent group that corresponds to
$w(\fn_+)$, then $U^w\cap U^{{\rm id}}$ acts freely on $U^{{\rm
id}}$, hence co-freely on $\CC[U^{{\rm id}}]$. Since $U^w\cap
U^{{\rm id}}$ and $U^{{\rm id}}$ are affine spaces, this translates
into the co-freeness of the action of the Lie algebra
$Lw(\fn_+)\cap\hat{\fn}_+$ on $\CC[U^{{\rm id}}]$. The passage to
the space of  loops in $U^{{\rm id}}$ is straightforward.

$\WW^w_{\nu_0/z}$ being co-free as an
$Lw(\fn_+)\cap\hat{\fn}_+$-module means that $(\WW^w_{\nu_0/z})^c$
is a free module over the ``opposite'' subalgebra,
$Lw(\fn_-)\cap\hat{\fn}_-$. Therefore,
\begin{equation}
\label{dual-wak-free-subsp}
(\WW^w_{\nu_0/z})^c=U(Lw(\fn_-)\cap\hat{\fn}_-)\otimes_{\CC}U
\end{equation}
for some graded subspace $U\subset (\WW^w_{\nu_0/z})^c$. The space
of co-invariants is
\[
\left((\WW^w_{\nu_0/z})^c\right)_{Lw(\fn_-)\cap\hat{\fn}_-}=U.
\]
Dualizing back one obtains the space of invariants
\[
(\WW^w_{\nu_0/z})^{Lw(\fn_+)\cap\hat{\fn}_+}=U^c
\]
or, in terms of formal characters,
\[
({\rm ch}\WW^w_{\nu_0/z})^{Lw(\fn_+)\cap\hat{\fn}_+}={\rm ch}U^c=
\frac{{\rm ch}\WW^w_{\nu_0/z}}{{\rm ch}U(Lw(\fn_-)\cap\hat{\fn}_-)}.
\]
One has the character formula
\begin{align}
&{\rm ch}U(Lw(\fn_-)\cap\hat{\fn}_-)=\nonumber\\
&\prod_{\begin{array}{c}\alpha\in
w(\Delta_+)\cap\Delta_+\\
n\geq 0 \end{array}}(1-e^{-\alpha-n\delta})^{-1}
\prod_{\begin{array}{c}\alpha\in
w(\Delta_+)\cap\Delta_-\\
n> 0 \end{array}}(1-e^{-\alpha-n\delta})^{-1}\nonumber
\end{align}
Dividing (\ref{char-of-wak-mod-restr}) by the latter we obtain
\begin{align}\label{char-of-space-almost-sing-vect}
&({\rm ch}\WW^w_{\nu_0/z})^{Lw(\fn_+)\cap\hat{\fn}_+}=\\
&e^{w\circ\nu_0}\prod_{\begin{array}{c}\alpha\in
w(\Delta_-)\cap\Delta_+\\
n\geq 0 \end{array}}(1-e^{-\alpha-n\delta})^{-1}
\prod_{\begin{array}{c}\alpha\in
w(\Delta_-)\cap\Delta_-\\
n> 0 \end{array}}(1-e^{-\alpha-n\delta})^{-1}\nonumber
\end{align}
Since $(\WW^w_{\nu_0/z})^{\hat{\fn}_+}\subset
\WW^w_{\nu_0/z})^{Lw(\fn_+)\cap\hat{\fn}_+}$, it follows that if
$v_{\mu}\in (\WW^w_{\nu_0/z})^{\hat{\fn}_+}$ is a singular vector of
weight $\mu$, then $\mu\in w\circ\nu_0+w(Q_+)$. Equivalently
\[
\mu=w\circ(\nu_0+\alpha)\text{ for some }\alpha\in Q_+.
\]
On the other hand, the block decomposition (\ref{block-decomp})
implies
\[
\mu\in W\circ\nu_0.
\]
Therefore $\nu_0+\alpha\in W\circ\nu_0$, but $\nu_0$ being dominant,
this may  happen only if $\alpha\in Q_-$, which requires that
$\alpha\in Q_+\cap Q_-=\{0\}$, hence $\mu=w\circ\nu_0$. A glance at
(\ref{char-of-space-almost-sing-vect}) shows that in
$(\WW^w_{\nu_0/z})^{Lw(\fn_+)\cap\hat{\fn}_+}$ there is only one
vector of weight $w\circ\nu_0$, ${\bf 1}$; (\ref{what we
want-lemma-dual-does-to-wak}) and Lemma~\ref{what-dual-does-to-wak}
follow. $\qed$

In order to conclude our proof of
Lemma~\ref{another-spectr-seq-Kpq}, hence of
Theorem~\ref{main-theorem-intro}(2) we need to consider an arbitrary
character $\nu(z)$.

\subsection{ Proof of Theorem~\ref{main-theorem-intro}(2):
the case of an arbitrary character.} Let us again denote by $\cC^j$
the direct sum $\oplus_{w\in W^{(j)}}\WW_{\nu(z)}$. According to the
key Lemma~\ref{reform-wak-modul-proof}, we have to compute the
cohomology of the complex
\begin{equation}
\label{key-complex-again}
0\rightarrow\cC^0\rightarrow\cC^1\rightarrow\cdots\rightarrow\cC^{{\rm
dim}X}\rightarrow 0.
\end{equation}

For an arbitrary $\nu(z)$, each $\cC^j$ carries an increasing
conformal filtration $\{F^i(\cC^j),\;n\geq 0\}$. A spectral sequence
arises (yet another one!), $\{(E^{ij}_r,d_r)\}$, with
\[
E^{ij}_0=F^{-i}\cC^{i+j}.
\]
 By definition,
$(\oplus_{ij}E^{ij}_0,d_0)$ is precisely the complex of
Lemma~\ref{reform-wak-modul-proof} with $\nu(z)$ replaced with
$\nu_0/z$. This places is in the situation of a homogeneous
character, the case we have just considered. Hence
$\oplus_{ij}E^{ij}_1$ is as asserted by
Theorem~\ref{main-theorem-intro}(2). The higher differentials
vanish. Indeed, all higher differentials are morphisms of graded
$\ghat$-modules. The modules in question are direct sums of
irreducibles, and morphisms among them are determined by the images
of highest weight vectors. But on those the differentials vanish
according to the following dimensional argument: by construction,
\[
d_r(E^{ij}_r)\subset E^{i+r,j-r+1}_r,
\]
thus {\it decreasing} the conformal weight. On the other hand, the
class of the highest weight vector $[{\bf 1}_{w\circ\nu_0}]$ may be
mapped only at a linear combination of classes  $[{\bf
1}_{u\circ\nu_0}]$ with $u>w$ (another use of the construction of
the Cousin complex (\ref{cousin}), and $\nu_0$ being regular
dominant, the conformal weights of the latter are strictly {\it
greater} than that of the former.

That the spectral sequence converges follows easily from the fact
that it lies inside a finite width band $\{(i,j),0<i+j<{\rm
dim}X,i<0\}$. $\qed$

\subsection{Proof of Corollary~\ref{char-form-cor-intro}}
Since by definition
\[
\chi(\cL_{\nu(z)}^{ch})=\sum_{i=0}^{{\rm dim}X}(-1)^i{\rm ch}
H^i(X,\cL_{\nu(z)}^{ch}),
\]
Theorem~\ref{main-theorem-intro} (2)  gives
\[
\chi(\cL_{\nu(z)}^{ch}) = \sum_{w\in W}(-1)^{l(w)}e^{-\langle\nu_0-
w\circ\nu_0,\rho^{\vee}\rangle\delta}{\rm ch}\VV_{\nu(z)},
\]

The Euler character, $\chi(\cL_{\nu(z)}^{ch})$, is known
(Theorem~\ref{main-theorem-intro} (1))
\[
\chi(\cL_{\nu(z)}^{ch})=\sum\limits_{w\in W}(-1)^{\ell(w)}e^{w\circ
\nu_0}\times \prod_{\alpha \in
       \widehat{\Delta}^{re}_+}
(1-e^{-\alpha})^{-1}.
\]
Therefore
\[
{\rm ch}\VV_{\nu(z)}=\sum\limits_{w\in W}(-1)^{\ell(w)}e^{w\circ
\nu_0}\times \prod_{\alpha \in
       \widehat{\Delta}^{re}_+}
(1-e^{-\alpha})^{-1}\times\left(\sum_{w\in
W}(-1)^{l(w)}e^{-\langle\nu_0-
w\circ\nu_0,\rho^{\vee}\rangle\delta}\right)^{-1}.
\]
The rightmost factor can be further factored out as follows
\[
\sum_{w\in W}(-1)^{l(w)}e^{-\langle\nu_0-
w\circ\nu_0,\rho^{\vee}\rangle\delta}=\prod_{\alpha\in
\Delta_+}(1-e^{-\langle \nu_0+\rho,\alpha^{\vee}\rangle\delta});
\]
this is a well-known identity, cf. \cite{FH} p.399. Plugging the latter identity into the above
expression for ${\rm ch}\VV_{\nu(z)}$ we obtain the desired result. $\qed$

\bigskip

\footnotesize{T.A.:  RIMS,
Kyoto University, Kyoto 606-8502 JAPAN.
Email:    arakawa@kurims.kyoto-u.ac.jp

F.M.: Department of Mathematics, University of Southern California, Los Angeles, CA 90089, USA.
E-mail:fmalikov@usc.edu}


\begin{thebibliography}{99}
\bibitem[A1]{A2} T.~Arakawa, Vanishing of cohomology associated to
quantized Drinfeld-Sokolov reduction, Int.Math.Res.Not., no.15
pp.730-767, 2004.

\bibitem[A2]{A4} T.~Arakawa, Representation Theory of Superconformal Algebras and
the Kac-Roan-Wakimoto Conjecture, Duke Math. J., Vol. 130 (2005), No. 3, 435-478.

\bibitem[A3]{A} T.~Arakawa, Representation theory of $W$-algebras,
Invent.Math., {\bf 169} (2007), no.3, pp.435-478.

\bibitem[A4]{A3} T.~Arakawa, Characters of representations of affine Kac-Moody Lie algebras at
 the critical level, arXiv:0706.1817v2.

\bibitem[A5]{A5} T.~Arakawa, Associated varieties of modules over Kac-Moody algebras and $C_2$-cofiniteness of W-algebras, arXiv:1004.1554v1.


\bibitem[AChM]{AChM} T.~Arakawa, D.~Chebotarov, F.~Malikov,
Algebras of twisted chiral differential operators and affine localization of $\fg$-modules, to appear in {\em Sel. Math.},
arXiv:0810.4964

\bibitem[AG]{AG} S.~Arkhipov, D.~Gaitsgory, 
    Differential operators and the loop group via chiral algebras , posted on
arXiv:math/0009007 

\bibitem[BB1]{BB1} A.~Beilinson, J.~Bernstein. Localisation de
$\frak{g}$-modules. (French)  C. R. Acad. Sci. Paris Se'r. I Math.
292  (1981),  no. 1, 15--18.


\bibitem[BB2]{BB2} A.~Beilinson, J.~Bernstein.
 A proof of Jantzen conjectures. I. M. Gelfand Seminar, 1--50,
 Adv. Soviet Math., 16, Part 1, Amer. Math. Soc., Providence, RI, 1993.


\bibitem[BD1]{BD1} A.~Beilinson, V.~Drinfeld. Chiral algebras. American
Mathematical Society Colloquium Publications, 51. American
Mathematical Society, Providence, RI, 2004. vi+375 pp. ISBN:
0-8218-3528-9


\bibitem[BD2]{BD2} A.~Beilinson, V.~Drinfeld. Quantization of Hitchin's
integrable system and Hecke eigensheaves, preprint.
%available at http://www.math.uchicago.edu/~mitya/langlands.html.

\bibitem[BL]{BL} L.~Borisov, A.~Libgober,  Elliptic genera of toric varieties
 and applications to mirror symmetry.
Invent. Math. 140 (2000), no. 2, 453--485.


\bibitem[Feig]{Feig} B.~Feigin,  Semi-infinite homology of Lie, Kac-Moody and
Virasoro algebras (Russian), {\it Uspekhi Mat. Nauk} {\bf  39 }
(1984),  no. 2, 155-156;



\bibitem[FF1]{FF1} B.~Feigin, E.~Frenkel, Representations of affine Kac-Moody
algebras and bosonization, in: V.Knizhnik Memorial Volume,  L.Brink,
D.Friedan, A.M.Polyakov (Eds.), 271-316, World Scientific,
Singapore, 1990

\bibitem[FF2]{FF2} B.~Feigin, E.~Frenkel, Affine Kac-Moody algebras at the
critical level and Gelfand-Dikii algebras, in: Infinite Analysis,
eds. A.Tsuchiya, T.Eguchi, M.Jimbo, {\it Adv. Series in Math. Phys.}
{\bf 16} 197-215, Singapore, World Scientific, 1992

\bibitem[F1]{F1} E.~Frenkel, Wakimoto modules, opers and the center at the
critical level. Adv. Math. 195 (2005), no. 2, 297--404.


%\bibitem[F2]{F2} E.~Frenkel,
 %Lectures on the Langlands program and conformal field theory.
 %Frontiers in number theory, physics, and geometry. II, 387--533, Springer, Berlin,2007.


\bibitem[F2]{F3} E.~Frenkel, Langlands correspondence for loop groups,
Cambridge University Press, 2007

\bibitem[FBZ]{FBZ} Frenkel E., Ben-Zvi D.,  Vertex algebras and algebraic
curves. Second edition. Mathematical Surveys and Monographs, 88.
American Mathematical Society, Providence, RI, 2004;

\bibitem[FG1]{FG1} E.~Frenkel, D.~Gaitsgory, D-modules on the affine
Grassmannian and representations of affine Kac-Moody algebras, {\it
Duke Math. J.} {\bf 125} (2004) 279-327

\bibitem[FG2]{FG2} E.~Frenkel, D.~Gaitsgory,   Weyl modules and opers without
monodromy, arXiv:0706.3725


\bibitem[FG3]{FG3} E.~Frenkel, D.~Gaitsgory,  Local Geometric Langlands
Correspondence: the Spherical Case, posted on arXiv:0711.1132

\bibitem[FH]{FH} W~.Fulton, J.~Harris, Representation theory. A first
course, Springer-Verlag, New York 1991.

%\bibitem[FZ]{FZ} I.~Frenkel, Y.~Zhu, Vertex operator algebras associated to
%representations of affine and Virasoro algebras. Duke Math. J. 66
%(1992), no. 1, 123--168.

\bibitem[GG]{GG} W.L.~Gan, V.~Ginzburg, Quantization of Slodowy slices, {\em Int.Math.Res.Not.}, no.2, p.243-255 (2004).

\bibitem[GM]{GM} V.~Gorbounov, F.~Malikov,
 The chiral de Rham complex and positivity of the equivariant signatures of some loop spaces.
 Manuscripta Math. 113 (2004), no. 3, 359--370.

 \bibitem[GelMan]{GelMan} S.I.~Gelfand,  Yu.I.~Manin, Homological
 algebra, 2ed., Springer-Verlag, Berlin Heidelberg New York, 1999

\bibitem [GMS1]{GMS}  V.~Gorbounov, F.~Malikov, V.~Schechtman, Gerbes of chiral
differential operators. II. Vertex algebroids, Invent. Math. 155
(2004), no. 3, 605-680.

\bibitem[GMS2]{GMSII} V.~Gorbounov, F.~Malikov, V.~Schechtman,
 On chiral differential operators over homogeneous spaces.  Int. J. Math. Math. Sci.  26  (2001),  no.2, 83--106.

\bibitem[HBJ]{HBJ} F. Hirzebruch, Th. Berger, R. Jung, Manifolds and modular forms. With appendices by
Nils-Peter Skoruppa and by Paul Baum. Aspects of Mathematics, E20. Friedr. Vieweg $\&$ Sohn,
Braunschweig, 1992.

\bibitem[K]{K} V.~Kac, Vertex algebras for beginners. University Lecture
Series, 10. American Mathematical Society, Providence, RI, 1997.
viii+141 pp. ISBN: 0-8218-0643-2.

\bibitem[KK]{KK} V.G.~Kac, D.~Kazhdan, Sructure of representations
with highest weight of infinite dimensional Lie algebras, Adv.Math.
{\bf 34} (1979), pp.97-108




 \bibitem[Kash]{Kash}M.~Kashiwara,  Kazhdan-Lusztig conjecture for a symmetrizable Kac-Moody Lie algebra.
 The Grothendieck Festschrift, Vol. II, 407--433, Progr. Math., 87, Birkh\"auser Boston, Boston, MA, 1990.

 \bibitem[KashTan1]{KashTan1}M.~Kashiwara, T.~Tanisaki,
 Kazhdan-Lusztig conjecture for symmetrizable Kac-Moody Lie algebra. II.
   Intersection cohomologies of Schubert varieties.
   Operator algebras, unitary representations, enveloping algebras, and invariant theory
    (Paris, 1989), 159--195, Progr. Math., 92, Birkh\"auser Boston, Boston, MA, 1990.

 \bibitem[KashTan2]{KashTan2}M.~Kashiwara, T.~Tanisaki,
 Kazhdan-Lusztig conjecture for symmetrizable Kac-Moody Lie algebras. III.
    Positive rational case.
    Mikio Sato: a great Japanese mathematician of the twentieth century. Asian J. Math. 2 (1998), no. 4, 779--832.
%\bibitem[Kost]{Kost} B.~Kostant, On Whittaker vectors and representation theory,
%Inv.Math. {\bf 48} (1978) 101-184.

%\bibitem[Li]{Li} H.~Li, Vertex algebras and vertex Poisson algebras, Contemp.Math., {\bf 6} %(2004), no.1, p.61-110

\bibitem[Li2]{Li2}H.~Li, Abelianizing vertex algebras,
{\em Comm. Math. Phys.}{\bf 259} (2005), no. 2, 391?411. , arXiv:0409140

\bibitem[Lus]{Lus} G.~Lusztig, Hecke algebras and Jantzen's generic decomposition patterns, {\em Adv. in Math.},
{\bf 37} (1980), no.2, p.121-164

 \bibitem[M]{M} F.~Malikov, Verma modules over Kac-Moody algebras of rank 2,
Leningrad Math.J., vol.2 (1991), no.2, p.269-286.

%\bibitem[MFF]{MFF} F.G.~Malikov, B.L.~Feigin, D.B.~Fuchs, Singular
%vectors in Verma modules over Kac-Moody algebras, Funct.Anal. i ego
%Pril. {\bf 20} (1986) no.2, pp.25-37.

%\bibitem[MS]{MS} F.~Malikov, V.~Schechtman, Chiral de Rham complex II,
%{\em D.B.Fuchs' 60 Anniversary Volume}, Amer.Math.Soc.Transl., {\bf 194} (1999), 149-188

\bibitem[MSV]{MSV} F.~Malikov, V.~Schechtman, A.~Vaintrob,  Comm. in Math. Phys. {\bf 204} (1999), 439-473



\bibitem[RCW]{RCW} A.~Rocha-Caridi, N.~Wallach, Projective modules over graded Lie algebras.I,
Math.Z. {\bf 180} (1982), pp.151-177.

\bibitem[V1]{V} A.~Voronov, Semi-infinite homological algebra,
Invent.Math. {\bf 113} (1993) pp.103-146.

\bibitem[V2]{V2}A.~Voronov, Semi-infinite induction and Wakimoto modules, {\em Amer. J. Math.}{\bf 121} (1999), no. 5, 1079-1094, q-alg/9704020.

\bibitem[W]{W} M.~Wakimoto Fock representations of the affine Lie algebra
 $A\sp {(1)}\sb 1$. Comm. Math. Phys. {\bf 104} (1986), no. 4, 605--609.

\bibitem[Witt]{Witt} E.~Witten,
Two-dimensional models with $(0,2)$ supersymmetry: perturbative aspects.
Adv. Theor. Math. Phys. 11 (2007), no. 1, 1--63.


\bibitem[Zhu]{Zhu} Y.~Zhu, Modular invariance of characters of vertex operator
algebras. J. Amer. Math. Soc. 9 (1996), no. 1, 237--302.
\end{thebibliography}
\end{document}